\documentclass[10pt]{article}

\usepackage[utf8]{inputenc}
\usepackage[english]{babel}
\usepackage{hyperref}
\hypersetup{
    colorlinks=false,    
    pdfborderstyle={},   
    pdfborder={0 0 0},  
    linkcolor=black,
    citecolor=black,
    filecolor=black,
    urlcolor=black
}

\usepackage{amsmath, amssymb, bm}

\usepackage{graphicx}
\usepackage{booktabs}
\usepackage{multirow}
\usepackage{makecell}
\usepackage{enumitem}

\usepackage[margin=1.5in]{geometry}
\makeatletter
\renewcommand\footnoterule{%
  \vspace{0.7em}   
  \hrule width 0.4\columnwidth height 0.4pt
  \kern 0.3em      
}
\makeatother
\usepackage[hang]{footmisc}
\usepackage[utf8]{inputenc}
\usepackage[english]{babel}
\usepackage[backend=biber,style=apa]{biblatex} 
\usepackage{amsmath}
\usepackage{booktabs}
\usepackage{float}
\usepackage{makecell}
\usepackage{tabularx} 
\usepackage[backend=biber,style=apa]{biblatex}
\usepackage[colorinlistoftodos]{todonotes}

\begin{document}

\title{\textbf{Preference-Based Optimisation for Integrated Design and Group Decision-Making}}

\author{
A.R.M. Wolfert\thanks{Corresponding author: A.R.M. (Rogier) Wolfert: \texttt{a.r.m.wolfert@tudelft.nl}}\\[0.5em]  
\small Department of Algorithmics, Faculty of Mathematics and Computer Science\\[-0.3em]  
\small Delft University of Technology, The Netherlands
}

\date{}
\maketitle

{
\renewcommand{\baselinestretch}{0.96}\normalsize
\vspace*{-2em}
\begin{abstract}
\noindent
Conventional multi-objective optimisation approaches (e.g., MOO-CP or MIP) typically fail in group decision-making by aggregating heterogeneous objectives without a valid preference foundation, producing Pareto-dominance sets instead of a unique, actionable decision. Because only humans can define objectives, their preferences constitute the only legitimate basis for decision-making. Accordingly, four essential conditions for complex, multi-level design–decision systems are established: (1) \textbf{Preference-Key} — all objectives, constraints, and trade-offs are evaluated within a unified preference domain using valid preference function modelling (PFM); (2) \textbf{Integration} — feasible system performance (“object-capability”) and acceptable actor preferences (“subject-desirability”) coexist within a single design–decision space; (3) \textbf{Association} — actors can freely specify their individual preferences and weights, enabling consistent aggregation towards group-optimal decision-making; and (4) \textbf{Uniqueness} — the solver identifies a single best-fit solution with maximum aggregated preference.

\noindent
The open design and decision system (ODESYS) methodology, employing the \textbf{IMAP} (Integrative Maximisation of Aggregated Preferences) optimisation framework together with its Intergenerational Solver (\textbf{IGS}) family, satisfies these conditions and enables fully integrated multi-objective design optimisation and multi-criteria decision-making (MODO-MCDM). Its extension within the \textbf{ODESYS/FIVES} formulation, introduced here, broadens applicability to a wider range of constrained design–decision problems while providing conceptual clarity through the explicit FIVES structure. Moreover, it operationalises affine preference aggregation while preserving equivalence with previously validated ODESYS 1.0 results. By mapping all system behaviour into a unified, goal-oriented preference-performance domain and directly aggregating actor preferences via PFM, ODESYS/FIVES/IMAP-IGS confronts system complexity and delivers a single, best-fit solution—even for highly constrained combinatorial group decision problems. In doing so, it guarantees feasible and acceptable outcomes, providing pure group design–decision support that surpasses conventional methods, which may generate numerical results yet lack group decision validity.

\noindent
Two illustrative applications—a marine installation and a multi-constrained vessel allocation problem—demonstrate how ODESYS/FIVES produces a unique, integrative preference-performance-based solution that maximises aggregated stakeholder preferences across system levels. This approach effectively transforms multi-objective optimisation into pure group decision-making, achieving a best-fit-for-common-purpose through the synthesis of performance possibilities and stakeholder purpose within socio-physical reach.
\end{abstract}

\vspace*{1.0em}


\vspace{0.2em}

\noindent\textbf{Article history:} Compiled June~30,~2026
\normalsize
}


\section*{Introduction}

Across systems engineering, design, and decision science, multi-objective optimisation (MOO) is widely recognised as a core group decision-support activity for confronting complexity. Contemporary problems involve heterogeneous objectives—e.g., availability, affordability, sustainability, and aesthetic value—evaluated by multiple stakeholders under uncertainty and dynamic conditions. Optimisation must therefore move beyond purely physical performance metrics and explicitly incorporate stakeholder preferences (utility or value), as only these provide a valid basis for goal-oriented choice. Separation between optimisation and preference evaluation undermines the design–decision outcome, since only the human-subject is goal-oriented and can set objectives, not the physical-object.This aligns with analyses emphasising that optimisation and preference elicitation cannot be meaningfully separated in decision-making processes \parencite{french2023reflections}.\\  

Foundational work, such as Decision-Based Design (DBD) \parencite{Hazelrigg1998}, Value Engineering \parencite{King2000ValueEngineeringTheoryPractice}, and utility-based design formulations \parencite{Thurston2011}, highlighted the necessity of preference-aware MCDM. However, these approaches fail to embed preferences as a rigorous basis for group decision-making and rely on numerical optimisation without grounding decisions in meaningful preference structures \parencite{french2023reflections}. Preference measurement was formalised through Barzilai’s Preference Function Modelling (PFM) theory \parencite{Barzilai2010Foundations, Barzilai2022PureEconomics}, demonstrating that preferences reside in a one-dimensional affine space and admissible transformations are strictly affine. Despite this, many aggregation practices remain mathematically invalid or decision-ambiguous. Recently, \parencite{Wolfert2026UniquePreference} demonstrated mathematically valid aggregation of preferences in MCDM using PFM.\\  

Classical systems engineering design optimisation texts \parencite{MartinsNing2022, BlanchardFabrycky2021} and management science references \parencite{HillierLieberman2021} provide algorithmic multi-objective solutions (scalarisation, $\varepsilon$-constraint, Pareto heuristics), but they focus on numerical optima rather than preference‑based MCDM. Consequently, solutions often yield Pareto sets or single optima lacking contextual preference meaning. This echoes concerns that Pareto‑based outcomes, while mathematically neat, do not themselves constitute decision-valid solutions \parencite{french2023reflections}. Contemporary literature \parencite{Pajasmaa2025GroupDecision, Ferdous2024} highlights the need for structured frameworks that consistently integrates preference, acceptability, performance, and feasibility to identify a design–decision vector that maximises aggregated preference.\\  

In practice, constrained MOO problems are commonly addressed using mathematical programming formulations such as LP, MILP, and MINLP, as well as constraint programming (CP) and hybrid MOO–MCDM frameworks \parencite{MartinsNing2022, HillierLieberman2021}, typically operationalised via weighted additive models, $\varepsilon$-constraints, lexicographic ordering, Pareto set generation, or hierarchical pruning strategies.These methods explore performance spaces but do not deliver a unique, group decision-valid solution; they optimise measures, not human, goal-oriented preferences. Moreover, exploring the full feasible design space is inefficient since only the subset of designs aligned with stakeholder preferences is relevant for effective decision-making. Decision-making requires preferences, which are grounded in human purpose rather than derived solely from data or models \parencite{french2023reflections, Barzilai2022PureEconomics}.  \\  

Similarly, advances in parametric and computational design \parencite{Woodbury2010ElementsOfParametricDesign, Block2013DesignComputation} enable systematic exploration of large design spaces through structured parameters and constraint-based modelling. While effective for evaluating technical and performance criteria (structural performance, sustainability, constructability, etc.), these methods remain fundamentally performance-driven, do not internalise stakeholder desirability or acceptability, and cannot produce a unique, decision-valid best-fit solution. They generate candidate designs requiring external interpretation, highlighting the gap between exploration and actual decision-making.\\  

To substantiate these findings, the following sections overview contemporary MOO methods for group decision-making, highlight gaps in research and practice, and present a development statement, including four conditions for preference-based optimisation that enable pure group decision-making.\\

\textbf{Persistent separation of preferences and system performance} remains a key limitation of classical optimisation. Preferences are often introduced only after technical feasibility has been established (i.e., within classical a-posteriori MCDA practice), rather than being intrinsically embedded within a unified design–decision space, as commonly observed in systems engineering and performance-driven design frameworks \parencite{deWeck2011EngineeringSystems, BlanchardFabrycky2021}. This separation persists across several contemporary paradigms, including Bayesian optimisation \parencite{Ahmadianshalchi2023}, generative design systems \parencite{ChenXu2025}, and dynamic preference models \parencite{Arezoomand2021, Kim2014, Regenwetter2022, Saadi2024}. Although these approaches increasingly incorporate preference information, such preferences are typically introduced indirectly through surrogate modelling, scalarisation, post-processing, or iterative interaction, rather than being fundamentally integrated into the optimisation structure itself. Similar limitations appear in interactive evolutionary algorithms \parencite{Thiele2009, Branke2016} and Choquet-integral or constraint-reformulation approaches \parencite{Hou2020, CastellanosAlvarez2021}. Consequently, changes in feasibility, capability, desirability, acceptability, and evolving project conditions fail to propagate coherently throughout the decision process. Moreover, exploring the feasible design space alone is inherently inefficient, since only the subset aligned with stakeholder preferences is decision-relevant, leading to fragmented solution spaces potentially without a convergent best-fit outcome. Meaningful design–decision optimisation emerges only when subjective preferences and objective system performance are coherently unified within the optimisation process itself, thereby shifting from classical MOO with a-posteriori MCDA toward an integrated a-priori MODO (multi-objective design–decision optimisation) framework for MCDM \parencite{french2023reflections, Wolfert2023OpenDesignSystems}. 
\vspace{2.5em}

\textbf{A key missing element in MOO is genuine participation}, as stakeholder preferences are not preserved as independent decision structures within optimisation frameworks, but are instead transformed through aggregation, weighting, or consensus mechanisms that enforce commensurability. As a result, individual interests cannot be preserved in their full structural form within the optimisation process. Actor dynamics in real-world decision environments are therefore characterised by heterogeneous and often structurally conflicting preference orientations over shared system attributes, reflecting the socio-technical nature of design decisions and stakeholder involvement \parencite{Yang2022PriorGoverned, Zhilyaev2022, Arkesteijn2017}. Participative, negotiation-based, and group decision frameworks improve descriptive realism \parencite{DuJiao2022, Shavazipour2025, Chou2021, Perez2023, Qiao2024, Lagaros2023, Rahimi2022}, yet preferences are typically treated as episodic inputs rather than continuously evolving variables embedded within the optimisation process. As a result, optimisation and decision-making remain structurally separated in most existing formulations \parencite{AdekoyaHelbig2023, ChenXu2025}. While recent work \parencite{Zhilyaev2022, Wang2025} has taken steps toward integration, a fully unified representation of engineering performance, feasibility, stakeholder preferences, and acceptability within a single associative design–decision space remains largely absent from mainstream MOO practice. Consequently, the identification of a best-fit-for-common-purpose solution remains under-explored, with stakeholder needs only partially internalised and the notion of a free, best-for-project choice vulnerable to becoming curated and illusory. Stakeholders and actor dynamics in real-world decision environments are characterised by heterogeneous and structurally conflicting preference orientations over shared system attributes. Each stakeholder is, in principle, free to evaluate system attributes according to their own interests, incentives, and decision context \parencite{Wolfert2023OpenDesignSystems}, implying that preference formation is inherently individual and not required to be aligned across actors. However, many multi-objective optimisation (MOO) formulations implicitly assume a shared or commensurable preference structure, typically operationalised through common objectives or aggregated weighting schemes, thereby masking underlying heterogeneity. In practice, actors may not only differ in preference intensity but also in directional orientation with respect to the same attribute, where a single system variable may represent benefit, cost, or risk depending on the stakeholder perspective. Consequently, preference information is generally non-aligned and non-commensurable across actors, forming heterogeneous preference geometries that must be explicitly represented at the actor level prior to any aggregation or optimisation.
\vspace{0.5em}

\textbf{The absence of a unified preference domain} for heterogeneous objectives remains a core limitation in multi-objective optimisation (MOO). Existing methods often resort either to hierarchical or semi-multi-objective formulations, or to monetisation as a surrogate for commensurability, despite longstanding critiques from decision theory, cost–benefit-oriented design, and socio-technical research \parencite{HirschHadorn2022, DuJiao2022}. Value is not solely monetary. Yet monetisation does not model preference; rather, it removes preference structure from the decision space and replaces it with scalar arithmetic convenience. Because monetary aggregation is not invariant under the affine structure of preferences, it collapses contextual, relational, and qualitative meaning into a single scalar representation \parencite{Barzilai2022PureEconomics}. Moreover, this overlooks the fact that monetary values originate in human judgments about worth and “value-for-money”. Once embedded in a cost-only optimisation, however, they are treated as absolute, additive, and universally comparable metrics, such that the optimisation no longer represents preference structure but instead executes pre-encoded value assumptions. Other non-monetised methods, such as hybrid MOO–MCDM frameworks, typically rely on a posteriori MCDM steps in which heterogeneous objectives are reconciled through method-dependent procedures rather than through a unified preference integration \parencite{Ferdous2024, ChenXu2025}. This fragmentation has been widely recognised as a persistent limitation in the development of coherent MCDM methodologies \parencite{french2023reflections}. Although explicit preference modelling is increasingly recognised as essential, systematic use of continuous, individually weighted preference functions grounded in rigorous measurement theory—and applicable across economic, technical, environmental, and social dimensions—remains rare. Notable recent efforts \parencite{Wang2025} integrate surrogate-assisted search with preference cues for expensive MOO problems while avoiding ideal point estimation, yet they remain objective-anchored and set-oriented rather than providing a fully associative, decision-valid preference domain \parencite{Lee2011, Messac1996, Arezoomand2021, Saadi2024, Wang2025}. In practice, heterogeneous objectives are treated in isolated spaces, with commensurability enforced via scaling, weighting, or dominance filtering rather than a unified preference domain. This confirms that the lack of a common preference domain is a persistent, cross-methodological challenge rooted in fundamental decision-theoretic limits, not mere algorithmic detail \parencite{Gunantara2018}.
\vspace{0.5em}

\textbf{Unique preference aggregation} is central to MOO–MCDM, as its mathematical validity determines the meaningfulness of decision outcomes. Nevertheless, most classical and post-2010 MOO–MCDM frameworks rely on normalisation, weighted sums, surrogate scores, or composite indices \parencite{Ferdous2024, AdekoyaHelbig2023, Kaddani2017, Dehshiri2022, Zeng2025}, which from a measurement-theoretic perspective, implicitly assuming ratio- or interval-scale properties that preference data do not possess. Similar concerns have been raised about the methodological fragility of weighted-sum and scaling-based aggregation techniques \parencite{french2023reflections}. Preference, however, is not a physical property but a subjective construct of the mind, expressing an individual’s free ordering of alternatives within a given context and defining a decision space that is inherently relational, individual, and situation-dependent. This is consistent with the broader understanding that preferences are constructed, context-dependent judgements rather than objective quantitative measures \parencite{french2023reflections}. More expressive methods—such as discrete choice models, network-based decision frameworks, and probabilistic preference embeddings \parencite{Sha2023}—increase representational or situational richness without resolving this foundational measurement inconsistency. Likewise, weighted aggregation with partial preference information, fuzzy–stochastic methods, and Pareto-front transformations continue to rely on normalisation, fuzzy membership functions, or transformed dominance relations—approaches that are operationally convenient yet methodologically fragile. Consequently, mathematically admissible and preference-consistent aggregation remains rare. Recent analyses \parencite{Pajasmaa2025GroupDecision} confirm that none of these approaches simultaneously satisfy the axioms of rigorous preference function modelling (PFM), which establishes that only differences between preference values are meaningful, admissible transformations are affine, and aggregation must preserve zero-reference stability and commensurability \parencite{Barzilai2010Foundations, Barzilai2022PureEconomics, Wolfert2026UniquePreference}.
\vspace{0.5em}

\textbf{Pareto remains a dominant epistemic anchor} in multi-objective optimisation \parencite{Marler2004}, including contemporary preference-aware variants \parencite{Ahmadianshalchi2023, Zhao2025, Pajasmaa2025GroupDecision}. Reviews \parencite{Gunantara2018} confirm that Pareto-dominance-based approaches are key within MOO, with preferences typically introduced post hoc via weights, constraints, or surrogate indices. This reliance on Pareto-type sets aligns with broader concerns that MCDM frameworks frequently cannot deliver unique decision outcomes \parencite{french2023reflections}. Many MOO methods result in Pareto sets requiring a posteriori selection, negotiation, or heuristic filtering. Pareto optimality identifies non-inferior (admissible) alternatives but cannot produce a unique, actionable solution—despite real-world multi-criteria design and decision-making requiring a single best-fit-for-common-purpose outcome. In this sense, Pareto optimality does not constitute a decision criterion, as it defines a mathematical dominance-based ordering in objective space rather than a cardinal preference-based selection in decision space. Consequently, without explicit preference information or a well-defined decision objective, the selection problem remains underdetermined, and no actionable “what-to-choose” can be inferred from the Pareto set alone. Even evolutionary multi-objective algorithms such a the well-known NSGA-II or NSGA-III \parencite{Deb2002, Deb2014}, or advanced extensions—such as interactive preference-guided methods \parencite{Thiele2009}, heuristic Pareto-front management \parencite{Kesireddy2024}, and hybrid aggregation-transform techniques \parencite{Zeng2025}—still produce sets or ranked subsets. Moreover, as the number of objectives increases, Pareto dominance rapidly loses discriminative power, because large Pareto sets often become indistinguishable without explicit preference modelling and are therefore treated as equivalently acceptable, despite stakeholder settings requiring a single acceptable and feasible optimum \parencite{Saad2018, Bai2015, Golany2006}. Consequently, Pareto optimality functions primarily as an ordered filtering mechanism rather than a decision-making concept. This limitation is structural rather than methodological. Conflicting preference directions on shared performance objectives and non-monotonic preference functions expose the limitations of Pareto dominance in maintaining preference consistency under such conditions. As a result, Pareto optimality is a purely ordinal notion of non-dominance in objective space and does not encode preference structure. Overall, Pareto-based analysis is incomplete as a decision framework because it operates purely in ordinal objective space without preference aggregation. Decision-making instead requires an explicit, preference-consistent aggregation operator in preference space, where solutions are evaluated via cardinal preference structures rather than dominance relations. Within such a framework, preference-consistent aggregation operators can yield unique optima that may lie outside the classical Pareto front, precisely because they optimise preference structure rather than non-dominance.
\vspace{0.5em}

\textbf{Real-world preferences are rarely linear or monotonic} because practical systems operate under constraints, thresholds, and diminishing or reversing returns. More is not always better, nor is less always preferable. For example, project scheduling involves trade-offs between acceleration risk and delay cost; urban density yields efficiency gains up to congestion thresholds; investment and innovation generate returns until coordination and complexity costs dominate; and spatial planning value increases only up to infrastructure, regulatory, or saturation limits. These examples reflect a general principle: decision-making seeks context-dependent optima rather than unbounded numerical MOO based on implicit linear or monotonic preference assumptions \parencite{french2023reflections, Barzilai2022PureEconomics} . Classical MI(N)LP formulations typically rely on convexity, monotonicity, and boundary-optimality, ensuring solutions occur at extreme points of the feasible polytope. Under non-monotone or piecewise preference functions, these structural interpretations break down, even though such models remain representable via mixed-integer linearisation. However, exact solvers such as CPLEX and Gurobi guarantee global optimality only within well-defined representable classes (e.g., MILP, convex QP), not for arbitrary non-monotone preference structures unless explicitly reformulated. While heterogeneous preferences across decision-makers are well recognised, the implications of conflicting monotonicity on a shared objective for Pareto-based decision-making are not typically explicitly formalised in these optimisation frameworks.

Even in single-objective optimisation (SOO), the preference structure is typically left implicit, effectively embedding unverified assumptions of linearity or monotonicity \parencite{Hazelrigg1998, MartinsNing2022} . This limits the applicability of these methods in real-world decision contexts, often requiring ad hoc linearisation techniques to ensure solvability. Preferences therefore exhibit saturation, asymmetry, acceptable ranges, tipping points, and multiple regimes of desirability. Consequently, the assumption of monotonic utility underlying classical optimisation and Pareto frameworks fails to represent actual decision behaviour, as real-world value landscapes are inherently nonlinear, non-convex, and often non-monotonic. In such settings, heterogeneous stakeholders introduce conflicting directional preferences on shared performance objectives, further decoupling Pareto efficiency from perceived desirability. When considered in high-dimensional multi-stakeholder objective spaces, this leads to Pareto degeneracy, where a large proportion of solutions become non-dominated and dominance loses discriminative power, as widely observed in many-objective optimisation literature (e.g. \parencite{Deb2014}).
\vspace{0.5em}

\textbf{Emerging preference-based optimisation} research challenges traditional dominance-based constructs, including ideal and nadir points \parencite{Zhao2025}, signalling a shift toward more explicitly preference-aware optimisation paradigms. Advanced approaches—such as interactive evolutionary algorithms \parencite{Branke2016}, constraint-enhanced preference integration \parencite{Hou2020}, and expensive preference-guided MOO methods (PMEGO) that avoid ideal points \parencite{Wang2025}—more effectively incorporate stakeholder objectives and capture complex criterion interactions. However, even these state-of-the-art methods do not satisfy the axioms of mathematically valid preference aggregation \parencite{Pajasmaa2025GroupDecision, Wolfert2026UniquePreference}. They remain fundamentally Pareto-set- or interaction-based: preferences steer, restrict, or filter candidate solutions, but are not structurally inseparable from system behaviour and do not produce a unique, decision-theoretically grounded outcome. Despite growing recognition of preference awareness, pure preference-based optimisation—mathematically valid, preference-meaningful, non-monetary, non-Pareto-based, and capable of yielding a single best-fit-for-common-purpose solution—remains absent from mainstream design and project management practice. To date, and to the best of current knowledge, ODESYS/IMAP \parencite{Wolfert2023OpenDesignSystems, VanHeukelum2024} is the only design–decision framework that explicitly operationalises integrative and associative preference-based MODO through the a priori use of preference function modelling, including unique preference aggregation, which enables a consistent and well-defined MCDM evaluation. This paper further extends the approach to a broader range of group design–decision applications, opening it to increasingly complex systems and addressing higher levels of system complexity.

\subsubsection*{Research \& Practice Gaps}

Across research and practice, persistent limitations remain in multi-objective design optimisation for group decision-making. Conceptually, preference-based design is widely acknowledged as necessary, yet current methods are incomplete. Pareto-dominance approaches structure search but cannot yield unique, decision-valid outcomes; preferences often guide exploration without full integration; aggregation methods frequently violate measurement-theoretic principles; and decision-making is deferred to episodic heuristics rather than being continuously embedded. Frameworks such as ODESYS introduced a unified preference-over-performance domain, formally integrating individual and acceptable stakeholder preferences with feasible system performance. Nevertheless, it remains partially objective-anchored, and its application to highly constrained socio-physical systems is under-explored. No framework yet fully realises a dynamically adaptive, mathematically rigorous preference-driven design–decision space capable of integrating system capability, feasibility, stakeholder desirability, and acceptability into a single best-fit-for-common-purpose solution.

In parallel, industry and public-sector practice faces a similar challenge. Existing methodologies remain fragmented, privileging either technical potential ('what can') or stakeholder desire ('what is wanted') without resolving their integration. Object-driven optimisation models technical performance but cannot generate meaningful best-fit-for-common-purpose outcomes; subject-driven methods—such as a posteriori ranking of pre-generated alternatives—simulate choice yet cannot guarantee actionable, feasible outcomes. Stakeholders are often confined to curated options, limiting decision freedom and potentially compromising results. Practitioners therefore require methods that embed capability and desirability within a shared, associative decision space, enabling transparent, explicitly substantiated decision-making while confronting complex socio-technical constraints. This preserves individual design freedom alongside collective alignment and legitimacy,ensuring that the conditions for identifying a best-fit solution are established a priori.

Taken together, both research and practice gaps highlight the absence of frameworks capable of concurrent, associative, and fully integrative decision-making, in which stakeholder preferences, system feasibility, and performance converge to produce a unique, best-fit-for-common-purpose solution.

\subsubsection*{Methodological Foundation}

The limitations identified in both research and practice motivate the evolution of the early ODESYS -- open design \& decision systems -- MODO framework \parencite{Wolfert2023OpenDesignSystems, VanHeukelum2024}, together with the IMAP (Integrative Maximisation of Aggregated Preferences) optimisation method and its Intergenerational Solver (IGS). Accordingly, there is a need for an open, integrative design and group decision-making MODO methodology that enables complex systems development across all relevant system levels, grounded in mathematically pure preference-based optimisation over performance and preference dimensions, and supporting open, human-centered group design and decision-making.

A rigorous formulation requires that such a preference–performance MODO method yields a unique, best-fit-for-common-purpose MCDM outcome if and only if it satisfies the following \textbf{four formal conditions}, which together define integrated design and group decision-valid preference-based optimisation.

\subsubsection*{Condition 1 — Preference-Key}\vspace{-0.7em}
A MODO method must represent and evaluate all performance objectives, constraints, and trade-offs within a single unified preference domain using mathematically consistent preference function modelling. Such methods treat preference as the key for decision-making and yield a unique best-fit-for-common-purpose outcome corresponding to the integrative maximisation of aggregated preference. This is realised through the a priori use of preference function modelling, whereby all design- and decision-relevant measures (including performance objectives and/or constraints) are transformed into preferences before optimisation.

\subsubsection*{Condition 2 — Integration}\vspace{-0.7em}
A MODO method must integrate subject-preference with object-performance. Feasible system performance (capability: “what it can”) and acceptable stakeholder preferences (desirability: “what they want”) must reside within a single integrative design-decision solution space. The immediately given physical reality—representing the system’s degrees of freedom (“infinite” possibilities)—must be coherently unified with human free will, expressed through individual actors’ goal-oriented preferences (common purpose); only under this integration can multi-objective optimisation produce a consistent design–decision synthesis.

\subsubsection*{Condition 3 — Association}\vspace{-0.7em}
A MODO method must be participative so that each individual stakeholder decision-maker is free to specify their own preference representation, including their local weighting and their acceptability limits within the feasible performance space. This formulation enables coherent cooperation and active participation across multiple actors and performance dimensions, defining a best-for-group decision as the maximisation of aggregated individual preferences over that space. Only within such an associative framework can a group-optimal decision emerge: when each actor advances their legitimate interests yet is willing to concede on pure self-interest, the group outcome achieves a higher aggregated preference than any purely individualistic design-decision alternative.

\subsubsection*{Condition 4 — Uniqueness}\vspace{-0.7em}
A MODO method must yield a unique and consistent best-fit-for-common-purpose solution —i.e., a best-fit design–decision point— that is invariant under admissible transformations of preferences, ensuring that the integrative maximisation of aggregated preference remains preference-meaningful and therefore valid, explicit, and transparent for open design and group decision-support. \\

Taking these four conditions as a starting point and building on recent ODESYS/IMAP 1.0 developments, the novel ODESYS in FIVES structure elevates preference-guided optimisation into a fully associative, design–decision-valid method for complex systems development. While ODESYS/IMAP 1.0 has demonstrated robustness both in academic courses and in practical implementations, including MSc-level applications at TU Delft, these formulations remained partially anchored in conventional multi-objective thinking. The full potential of a \textbf{pure preference-performance approach}—eliminating the objective layer and enabling fully integrative design and group decision-valid MCDM for highly constrained system contexts—is realised in the extended ODESYS operator structured in FIVES. Moreover, ODESYS in FIVES introduces the unique preference aggregation mechanism, the \emph{a-fine-aggregator}, which formally consolidates stakeholder preference information into a single preference-consistent decision criterion. This extension preserves previously validated ODESYS outcomes \parencite{Wolfert2023OpenDesignSystems} while significantly extending the methodological scope toward highly constrained socio-physical design–decision problems.

Building on this foundation, ODESYS moves beyond forced compromises inherent in conventional multi-objective, constraint-based, and Pareto-based optimisation approaches toward a genuinely best-fit-for-common-purpose outcome. By internalising human preference, desirability, and acceptability directly within the feasible system performance space, ODESYS transforms parametric design exploration into a group decision-valid synthesis process embedded in reflective practice.

\newpage
\section*{1. Open design and decision system}
Let a complex, multifaceted socio-physical problem comprising system performance functions indexed by $i$ and actor preference functions indexed by $k$, and defined over the endogenous decision vector $\mathbf{x}$, the exogenous parameter vector $\mathbf{y}$, and the system state at time instant $t$, be formulated as an open design and decision system (ODESYS) synthesis operator:

\begin{equation}
\mathrm{OD}\Big(P_{k,i}(\mathbf{F}(\mathbf{x},\mathbf{y},t));\; w'_{k,i}\Big), 
\label{OD operator}
\end{equation}

Here, $\mathbf{F}=\big(f_1..f_I\big)$ denotes the vector of system performance functions, where each $f_i(\mathbf{x},\mathbf{y},t)$ denotes the $i$th performance function, and $I$ is the total number of performance functions. These functions describe the physical behaviour of the system, i.e. its capabilities, representing what the system \emph{can} perform. The vector $\mathbf{x}$ represents the 'controllable' design-decision vector, comprising both product and process variables that can be deliberately selected and adjusted by the decision maker. The vector $\mathbf{y}$ is the 'uncontrollable' parameter vector, containing contextual variables, including external product and process conditions that influence system behaviour but cannot be controlled. The time variable $t$ may represent an explicit dynamic state, a fixed decision horizon $t = T^*$, or be omitted in static formulations, yielding $\mathbf{F}(\mathbf{x},\mathbf{y})$.

In contrast to traditional optimisation approaches, these performance functions over the design--decision vector do not constitute objectives in themselves, as systems are not inherently goal-oriented. Goal-oriented decision-making is enabled only by transforming the performance functions into actor-specific preference functions $P_{k,i}(\cdot)$. A preference expresses the relative desirability, value, or utility of a design alternative or decision option, thereby capturing what an actor \emph{wants} from the system. Preference functions reflect the preferences of actor $k=1,\ldots,K$ (with $K$ denoting the total number of actors) with respect to performance function $f_i$, yielding a total of $K \cdot I$ preference functions. Depending on the actor's interest, a preference function may be monotonically increasing (preferring higher performance values), monotonically decreasing (preferring lower performance values), or any combination. In this sense, preference functions generalise the traditional notions of objective maximisation and minimisation in a unified preference--performance space. Preference functions are freely specified by the individual actors themselves and may therefore adopt any (non-)linear or (non-)monotonic form. These preference functions reflect real-world stakeholders' individual performance interests and may consequently produce conflicting directional preferences over the same performance dimension in preference space.

Formally, the preference mapping can be interpreted as follows. Preference is a subjective and contextual construct that is not an absolute physical property but represents an individual's relative choices, thereby defining the decision space. The preference function $P_{k,i}$ maps an actor-independent performance function vector $\mathbf{F}$ to a corresponding preference score, reflecting the goal-oriented evaluation of actor $k$, such that $P_{k,i}(\mathbf{F}) \in [0,100]$, where $P_{k,i}$ acts component-wise on the $i$th element of $\mathbf{F}$. The relatively worst (least preferred) and best (most preferred) performance levels are assigned preference scores of 0 and 100, respectively, representing an arbitrary reference scaling since preferences possess only relative meaning within a one-dimensional affine space, i.e., all preference values are defined only up to relative differences and not in absolute terms \parencite{Wolfert2026UniquePreference}.

In addition to directional preference, actors may express the relative importance of performance dimensions through preference weights. This accounts for both (i) the actor’s overall importance in the decision process, represented by the global weight $w_k$ (with $w_k = 1/K$ in the case of equivalent stakeholder importance), and (ii) the local importance weight \(w_{k,i}\) assigned by stakeholder \(k\) to performance function \(i\), such that the effective contribution weight for each preference function is given by \(w'_{k,i} = w_k \cdot w_{k,i}\), with \(\sum_k w_k = 1\), \(\sum_i w_{k,i} = 1 \ \forall k\), and consequently \(\sum_{k,i} w'_{k,i} = 1\). When an actor has no interest in a particular performance function, the corresponding local weight is set to zero.   

The ODESYS system is subject to hard feasibility constraints defined in performance space, including their design-decision-related implications, via the \emph{feasibility} operator
\begin{equation}
g_f(\mathbf{F}) := \mathbf{F} - \overline{\mathbf{F}} \le 0,
\end{equation}
where $\overline{\mathbf{F}}$ denotes the performance-bound vector defining the feasibility envelope, representing externally imposed limits. These constraints are non-negotiable and cannot be influenced by the acting decision maker.

In contrast, the decision maker can define acceptability limits on performance functions through soft constraints expressed entirely in the preference function space $P_{k,i}(\mathbf{F})$. Any implicit performance restrictions arise from the actor's choice of reference intervals used to construct the preference functions, rather than from explicit feasibility constraints imposed on the performance functions. For each performance function $f_i$, an individual reference interval  $f_i \in [\, f_i^{\mathrm{loc}} \, ; \, f_i^{\mathrm{upc}} \,] \subseteq [\, f_i^{\min} \, ; \, f_i^{\max} \,]$ is defined, where $f_i^{\min}$ and $f_i^{\max}$ denote the global extrema over the feasible system space $\mathcal{S}_f$. The endpoints $(f_i^{\mathrm{loc}}, f_i^{\mathrm{upc}})$ correspond to preference scores of 100 (relative-best) and 0 (relative-worst), respectively, while intermediate values are mapped through the actor-specific preference function $P_{k,i}(\mathbf{F}) \in [0,100]$. (Note: in the case of linear preference functions defined over the full feasible range, i.e., $f_i^{\mathrm{loc}} = f_i^{\min}$ and $f_i^{\mathrm{upc}} = f_i^{\max}$, 
$P_{k,i}(\mathbf{F}) = 100 \cdot (f_i - f_i^{\min})/(f_i^{\max} - f_i^{\min})$.)

Further, the ODESYS system can be subject to additional soft constraints defining acceptability in preference space via the component-wise \emph{acceptability} operator:
\begin{equation}
\label{eq:acceptability}
g_a\big(P_{k,i}(\mathbf{F})\big) := P_{k,i}(\mathbf{F}) - \overline{P}_{k,i} \ge 0,
\end{equation}
where $\overline{P}_{k,i}$ denotes the minimum acceptable preference threshold in preference space, representing the actor-specific acceptability boundary (typically 0), although stricter cut-off levels may be imposed to reflect higher selection standards.

Finally, the complete design--decision system solution space, encompassing both feasibility and acceptability constraints, is defined as

\begin{equation}
\mathcal{S}_{f,a} := 
\left\{ 
\mathbf{x} \;\middle|\; 
g_f(\mathbf{F}) \le 0 \;\wedge\; 
g_a(P_{k,i}(\mathbf{F})) \ge 0 \;\wedge\; 
f_i(\mathbf{x}) \in [f_i^{\mathrm{loc}},\,f_i^{\mathrm{upc}}]
\right\}.
\label{Sol Space}
\end{equation}

where $P_{k,i}(\mathbf{F}) \in [\overline{P}_{k,i},100]$ denotes the actor-specific preference function over performance dimension $f_i$ (i.e., the acceptability envelope), with negotiable acceptability defined through the actor-specific performance reference intervals $[f_i^{\mathrm{loc}},\,f_i^{\mathrm{upc}}]$, which define the domain over which the preference functions are constructed. This set defines the admissible design–decision space prior to optimisation over the design–decision vector $\mathbf{x}$, i.e., the subset of decision vectors that induce feasible (hard-constrained) and acceptable (soft-constrained) system behaviour through the system performance mapping $\mathbf{F}$ and the associated preference functions $P_{k,i}(\mathbf{F})$.\vspace*{0.6em}

\textsc{Notes} to Equations (1)-(4): \vspace*{0.1em}

(1) By distinguishing system performance functions $\mathbf{F}$ from stakeholder preference functions $P_{k,i}(\mathbf{F})$, system behaviour is embedded within a design--decision framework that links the decision vector $\mathbf{x}$ to an objectively quantifiable performance space (i.e.\ measurable system outputs derived from physical and operational variables), and subsequently to a subjectively evaluated preference space in which stakeholder-specific assessments are formed. This distinction reflects the fact that system performance functions are grounded in objective, non-goal-oriented system behaviour, while the subjective interest assigned to those performance functions may differ across goal-oriented stakeholders. For example, traffic nuisance is an objectively measurable performance characteristic determined by design variables such as road width and traffic speed, yet its evaluation may be perceived differently in terms of acceptability or disturbance. Similarly, costs are derived from objective design quantities such as material usage or time, but their assessment is ultimately expressed in subjective perceived value for money. Finally, aesthetic appearance is inherently more subjective, but can still be systematically linked to objective design variables such as material choice and architectural features. The final assessment is, however, determined in the preference space.\vspace*{0.1em}

(2) The feasibility operator is evaluated in performance space through $\mathbf{F}$, while being enforced on the underlying design-decision and parameter variables $(\mathbf{x},\mathbf{y})$, thereby allowing constraints to be formulated either in performance space or directly in design-decision-parameter space via the system performance mapping. These are captured by $g_f(\mathbf{x},\mathbf{y},t) \le 0$, encompassing domain, logical, activity, sequencing, path, and physical (capacity) constraints, and may include equality constraints $h_f(\mathbf{x},\mathbf{y},t) = 0$. Overall, the feasibility operator provides a concrete instantiation of hard constraints, which are non-negotiable and arise from (a) performance-envelope limits (e.g.\ safety levels, external regulatory limits, and environmental conditions); (b) decision-structure constraints (e.g.\ availability, activity and logical coupling); (c) design-physics constraints (e.g.\ mechanical or (hydro)dynamics laws and bearing capacity).
Soft constraints, in contrast, are negotiable and act on the preference functions $P_{k,i}(\mathbf{F})$, reflecting stakeholder interests, priorities, and acceptable performance levels. By defining an individual preference function with preference levels between the minimum $\overline{P}_{k,i}$ (typically 0) and the maximum preference level (typically 100), over an acceptability interval $[\,f_i^{\mathrm{loc}}, f_i^{\mathrm{upc}}\,]$, the decision maker implicitly bounds acceptable system performance while ensuring that all acceptability reasoning remains within preference space.

\subsection*{2. ODESYS structure in FIVES}

\noindent
In complex socio-physical problems, multiple stakeholders or actors express preferences over a set of system performance functions $\mathbf{F}(\mathbf{x},\mathbf{y},t)$, as formulated through the ODESYS synthesis operator (see Equation~\eqref{OD operator}). Each performance function describes a capability of the system, i.e., a behaviour that the system can perform within its physical limits. Each preference function $P_{k,i}$ captures the desirability assigned by stakeholder $k$ to performance function $f_i$, together with its associated importance weight $w'_{k,i}$. Since stakeholders may hold different and potentially conflicting interests, a design–decision methodology is required that can consistently integrate these preferences while respecting both non-negotiable feasibility and negotiable acceptability constraints.
Within the ODESYS/FIVES MODO framework, the first four elements—capability, feasibility, desirability, and acceptability—define the integrated design–decision structure. The fifth element, solvability, is provided through the Integrative Maximisation of Aggregated Preferences (IMAP) method: a preference-based constrained optimisation principle that aggregates stakeholder preferences into a single design–decision criterion to be maximised. The resulting optimisation identifies the best-fit-for-common-purpose design–decision solution within the hard feasibility constraints and soft acceptability limits, over the solution space $S_{f,a}$ (see Equation~\eqref{Sol Space}).\\ 
Formally, the IMAP–MODO preference-based multi-objective design optimisation method for multi-criteria group decision-making is defined as:

\begin{equation}
\max_{\mathbf{x} \in \mathcal{S}_{f,a}} Z
= \mathbf{A}\big(z_{k,i}(\mathbf{x});\; w'_{k,i}\big),
\qquad
z_{k,i}(\mathbf{x}) := \mathcal{N}\!\big(P_{k,i}(\mathbf{F})\big)
\label{maxZ}
\end{equation}

where $\mathbf{A}$ is the \textbf{a-fine-aggregator}, a linear weighted centroid operator that aggregates actors’ normalised preference scores over their associated performance evaluations. The operator acts on the $z$-normalised preference scores $z_{k,i}(\mathbf{x})$, weighted by $w'_{k,i}$, inducing a preference-consistent aggregated ranking over the design–decision space, see \parencite{Wolfert2026UniquePreference} for details on unique preference aggregation using Preference Function Modelling (PFM)). The best-fit-for-common-purpose solution is defined as the unique design–decision vector $\mathbf{x}^*$ that maximises the aggregated preference operator:

\begin{equation}
\mathbf{x}^* = \arg\max_{\mathbf{x} \in \mathcal{S}_{f,a}} Z,
\qquad
Z = \sum_{k,i} w'_{k,i} \, z_{k,i}(\mathbf{x}),
\qquad
\sum_{k,i} w'_{k,i} = 1
\label{xStarr}
\end{equation}

Solvability is defined as the existence and computability of a design--decision vector $\mathbf{x}^*$ that maximises the aggregated preference $Z$, subject to hard feasibility constraints and soft acceptability limits; this constitutes the Integrative Maximisation of Aggregated Preferences (IMAP) optimisation method. Unlike classical notions of decidability, solvability in this context presupposes the existence of a computational decision-support mechanism and focuses on resolving a complex design--decision space rather than establishing theoretical computability alone within a classical objective-space formulation.

The resulting design--decision engine is the so-called Preferendus, a software system of preference-based optimisation solvers that identifies the optimal design--decision vector from Equation~\eqref{xStarr} (see the Data Availability statement and Preferendus applications in \parencite{Wolfert2023OpenDesignSystems,Teuber2025Odycon}). Note that when $\mathbf{x}$ is empty and all variables are exogenous, the MODO--IMAP optimisation formulation reduces to a pure MCDA (multi-criteria decision analysis) problem, where ranking is performed using the a-fine-aggregator.

Within Preferendus' IMAP method, preference-performance optimisation — the search for a best-fit — is operationalised through an Inter-Generational Solver (IGS), in which population-based metaheuristic search processes serve as exploration mechanisms within the constrained search space. In particular, genetic algorithms (GAs) and swarm intelligence methods belong to the broader class of population-based metaheuristics that provide alternative instantiations of the search component. More generally, the IGS architecture decouples stochastic search from preference-based decision aggregation by embedding interchangeable population-based metaheuristics within a higher-level preference-consistent aggregation layer. The current Preferendus implementation includes multiple IGS instantiations and is available as an open-source package, accessible via the links provided in the Data Availability section (e.g., BRKGA-based IGS instantiations for combinatorial and highly constrained optimisation problems, see \parencite{Timp2026}, or swarm-intelligence (SI)-based IGS instantiations for routing and logistics applications, see \parencite{Ozer2026}).

IMAP-IGS operates within a normalised preference space and handles contextual preference scores derived from the underlying preference functions. The resulting ranking dynamics are operationalised through the PFM-consistent aggregation operator ($Z$), which induces a design–decision representation enabling unique ranking and selection. The IGS architecture is characterised by three core components: (i) generation-to-generation preference propagation, (ii) a memory structure that preserves and updates preference-relevant state information across iterations, and (iii) a PFM-consistent a-fine-aggregator that maps population-level information into a unique preference-consistent decision within the constrained search space. IGS represents a family of algorithms that embeds IMAP scoring into the evolutionary fitness evaluation of population-based search methods. By replacing the standard scalar fitness function with the IMAP a-fine-aggregator, the framework preserves the underlying evolutionary operators while steering the search toward the group-preference optimum. This makes the approach broadly applicable to population-based metaheuristics, as it allows existing optimisation algorithms to be augmented with preference-consistent decision logic without modifying their core search dynamics.

In contrast to classical population-based algorithms, the IGS-GA family does not rely solely on generational population replacement for decision-making. Instead, it preserves the highest-ranked solution from the latest generation, ensuring at minimum a one-step intergenerational memory from $G_{n-1}$, within an intergenerational reference framework, enabling preference-consistent comparison under contextual and normalised preference rankings. This may be interpreted as a form of preference-based elitism within the search algorithm. Unlike classical elitism, which preserves solutions based solely on objective fitness, the IGS-GA preserves the highest-ranked solution according to the contextualised and normalised preference ranking. Extended implementations may additionally retain historical best-ranked solutions across multiple generations ($G_0$ to $G_{n-1}$) to further stabilise convergence and improve robustness against contextual rank fluctuations \parencite{VanHeukelum2024}.

Finally, IGS handles constraint violation in a preference-based manner. To guide the search when the evaluation pool consists entirely of infeasible or unacceptable solutions, the framework introduces a preference-consistent constraint-handling mechanism in which constraint violation is formulated as an explicit preference function. Feasible solutions receive a maximum relative preference score of 100, while infeasible solutions are assigned a monotonically decreasing preference value as a function of their constraint violation magnitude. Hard constraint violations are mapped to the lowest preference level in the preference scale, ensuring strict prioritisation of feasibility. Consequently, even within such populations, the aggregation operator can still discriminate and prioritise candidate solutions that are closer to the feasible–acceptable region. In this way, the IGS family explores the feasible–acceptable solution space $\mathcal{S}_{f,a}$ through aggregation of normalised preference scores, ultimately yielding the best-fit-for-common-purpose design-decision vector $\mathbf{x}^*$, i.e. the solution that maximises the aggregated preference score.\\

\vspace{-0.8em}
\subsubsection*{Summary --- ODESYS structure in FIVES}\vspace{-0.5em}
The ODESYS operator, formulated in Equation~\eqref{OD operator}, describes a multifaceted socio-physical problem comprising system performance functions indexed by $i$ and actor preference functions indexed by $k$, defined over the controllable design--decision vector $\mathbf{x}$ (endogenous), the uncontrollable parameter vector $\mathbf{y}$ (exogenous), and the system state at evaluation instant $t$. Within FIVES, this operator is structured as an integrated design and group decision-making MODO framework, implemented through the IMAP preference-based optimisation method and its IGS intergenerational solver, using the full set of equations defined in Sections~(1--2):

\begin{equation}
\label{eq:odesys_fives_summary}
\begin{aligned}
(1)\;& \text{Capability (performances)}      && \mathbf{F} \\
(2)\;& \text{Feasibility (hard constraints)} && g_f(\mathbf{F}) \le 0 \\
(3)\;& \text{Desirability (preferences)}     && P_{k,i}\big(\mathbf{F}\big) \\
(4)\;& \text{Acceptability (soft constraints)} 
&& g_a\big(P_{k,i}(\mathbf{F})\big) \ge 0  \\
(5)\;& \text{Solvability (IMAP-IGS)}
&& \displaystyle \max_{\mathbf{x}} \; Z
   = \mathbf{A}\big(z_{k,i}(\mathbf{x});\; w'_{k,i}\big)
\end{aligned}
\end{equation}

\vspace{0.3em}
\vspace{0.3em}
\textsc{Notes} to Equation \eqref{eq:odesys_fives_summary}: \vspace*{0.1em}

(1) The ODESYS/FIVES structure extends the ODESYS/IMAP 1.0 MODO framework to a broader class of group design–decision applications addressing socio-physical system complexity. Building on the validated robustness of predominantly design-driven ODESYS/IMAP 1.0, it preserves established outcomes while enabling an integrated representation of capability, feasibility, desirability, acceptability, and explicit IGS-based solvability across multi-actor and highly constrained problem settings.\vspace*{0.1em}

(2) In the ODESYS/FIVES formulation, objective functions are not treated as independent optimisation targets. Instead, goal orientation is expressed through stakeholder preference functions \( P_{k,i} \), which constitute the a priori MODO basis for design–decision-making. Best-fit solutions are obtained by maximising aggregated preferences rather than by direct optimisation of extrema in performance space or by numerical scaling or monetisation of objectives. This establishes a preference-based decision space as the primary domain for evaluation and optimisation, accommodating heterogeneous stakeholder interests and actor-specific weight structures within a consistent aggregation framework.\vspace*{0.1em}

(3) In the ODESYS/FIVES formulation, acceptability is defined in the preference domain via \( g_a(P_{k,i}(\mathbf{F})) \), but may also be implicitly embedded in the construction of the preference functions through actor-specific reference intervals \( [f_i^{\mathrm{loc}}, f_i^{\mathrm{upc}}] \). In this way, acceptable performance ranges can be incorporated directly into the preference structure, without requiring explicit constraint formulations. This dual representation enables flexible modelling of negotiable actor-specific limits alongside non-negotiable feasibility constraints.\vspace*{0.1em}

(4) ODESYS departs from classical multi-objective optimisation (MOO) paradigms by replacing set-valued optimality (e.g., Pareto fronts) with a single, preference-consistent design–decision solution \( \mathbf{x}^* \). Within the IMAP framework, this solution maximises aggregated preferences over the integrated performance–preference space, thereby providing a unified design–decision outcome rather than a set of non-dominated alternatives.\vspace*{0.1em}

(5) ODESYS can also be applied to single-objective problems. If a single preference weight \( w'_{k,i} = 1 \) and all others are zero, the IMAP aggregation reduces to a single-objective case in which the preference directly represents the underlying performance measure, yielding a solution analogous to classical \texttt{MINIMIZE} or \texttt{MAXIMIZE} formulations. However, in the presence of multiple stakeholders or conflicting preferences, even for a single performance dimension, classical single-objective formulations are generally insufficient to produce design–decision-valid outcomes.

\newpage
\vspace{0.6em}
\section*{3. Demonstrative ODESYS/FIVES problems}

In this section, we present two demonstrators to illustrate the ODESYS/FIVES framework (see Equations~\eqref{OD operator} and~\eqref{eq:odesys_fives_summary}). Demonstrator~\#1 concerns an integrative design--decision problem for the installation of a floating wind farm, which has been used as a demonstrator in several previous publications \parencite{VanHeukelum2023,VanHeukelum2024,Teuber2025Odycon}, as well as in the ODESYS book example presented in Chapter~8.5 \parencite{Wolfert2023OpenDesignSystems}, and through the links provided in the Data Availability section.

Demonstrator~\#2 concerns a dynamic, highly constrained combinatorial fleet-allocation decision problem for a marine service contractor, which led to the development of the ODESYS/FIVES-based Boskalis AlloDyn decision-support system. In this case, the methodology is likewise presented together with a modelling framework for preference-based group decision-making. The focus of this demonstrator is not on the operational outcomes produced by AlloDyn, but on formally describing the underlying multifaceted decision problem and its multiply constrained solution space (for corresponding results, including a comparative MOO benchmark, see the HVASP problem presented in \parencite{Timp2026}). As with Demonstrator~\#1, this paper does not discuss specific problem results. Instead, it describes the CP-IMAP preference-based optimisation method and introduces a proprietary IGS-GA instantiation, with corresponding IMAP-IGS implementation links provided in the Data Availability section.

The purpose of both demonstrators is to exemplify the generalised ODESYS/FIVES framework and its solution methodology, including preference-consistent aggregation, rather than to present detailed problem-specific results. The reader is referred to the aforementioned publications for the underlying design--decision results. While the results remain unchanged, the formulation presented here is novel and preserves the previously validated outcomes.

\subsubsection*{\textsc{Demonstrator}~\#1: Floating Wind Farm Problem}
\vspace*{-0.6em}
This demonstrator addresses a complex, multi-actor installation planning problem for a floating wind farm, in which engineering design choices and managerial decisions are tightly coupled with four performance functions of interest: project duration (make-span), total cost, fleet utilisation, and vessel emissions. The installation problem involves multiple vessels and two concurrent actors—the energy service provider and the marine contractor—with potentially conflicting interests, and is formulated as the ODESYS synthesis operator in Equation~\eqref{OD operator}, where $k = 1..2$ and $i = 1..4$. This formulation constitutes an \emph{integrative design--decision problem}, as it simultaneously addresses decision variables related to process management (e.g.\ vessel selection and deployment) and engineering design choices (e.g.\ anchor dimensions), involving two stakeholders, their associated preference functions, and four system performance functions. Moreover, in static formulations, time represents the evaluation instant of a fixed ODESYS system state and is therefore omitted, yielding $\mathbf{F}(\mathbf{x},\mathbf{y})$.\vspace*{0.6em}

The ODESYS structure in FIVES is defined as follows.

\paragraph{(1) Capability}
The system capability is represented by a four-dimensional system performance vector
\[
\mathbf{F}=
\big[
f_1(\mathbf{x}),
f_2(\mathbf{x},\mathbf{y}),
f_3(\mathbf{x}),
f_4(\mathbf{x},\mathbf{y})
\big],
\]
where the performance functions of interest respectively capture project duration (make-span) $f_1(\mathbf{x})$, installation cost $f_2(\mathbf{x},\mathbf{y})$, fleet utilisation $f_3(\mathbf{x})$, and vessel-related emissions $f_4(\mathbf{x},\mathbf{y})$. The project duration function $f_1$ is evaluated using a discrete-event simulation (DES), which determines the total installation time based on vessel types, installation rates, deck capacities, and anchor reloading operations. The cost function $f_2$ combines anchor manufacturing costs (Capex) with time-dependent vessel operating costs (Opex), where vessel operating times are obtained from the DES. The fleet utilisation function $f_3$ captures the strategic suitability of the selected vessels by quantifying the combined opportunity cost associated with deploying the selected vessels to the project rather than elsewhere. The vessel emissions function $f_4$ accounts for total project-related emissions by aggregating vessel-specific emission rates over their operational durations, as determined by the DES. The DES provides vessel operating times $t_i(\mathbf{x})$, which are used in $f_1$, $f_2$, and $f_4$. Temporal dynamics such as sailing times, installation sequencing, and vessel utilisation are therefore embedded implicitly within the performance functions, in particular in $f_1(\mathbf{x})$, $f_2(\mathbf{x},\mathbf{y})$, and $f_4(\mathbf{x},\mathbf{y})$.

The controllable design--decision vector is $\mathbf{x} = (x_1,\dots,x_5)$, while the uncontrollable parameter vector is $\mathbf{y} = (y_1,\dots,y_7,\mathbf{y}^{R},\mathbf{y}^{E})$. The full analytical expressions for $f_1$--$f_4$, including their dependence on $\mathbf{x}$ and $\mathbf{y}$, are provided in the Appendix.

\paragraph{(2) Feasibility}
System feasibility is enforced through a set of constraints that ensure:
- domain validity (variable bounds on the decision vector $\mathbf{x}$ and exogenous parameter vector $\mathbf{y}$),
- logical validity (feasible vessel selection), 
- and physical feasibility (the load-bearing capacity of the cable-anchor system). 
The feasible system solution space prior to optimisation over the design–decision vector $\mathbf{x}$ is defined as:
\[
\mathcal{S}_f
=
\left\{
\mathbf{x}
\;\middle|\;
g_f^{(m)} \le 0,\;
m = 0..2
\right\},
\]
Here, $g_f^{(0)}$ represents the domain constraints for the vectors $\mathbf{x}$ and $\mathbf{y}$, while constraints $g_f^{(1)}$–$g_f^{(2)}$ capture logical and physical constraints. All individual constraints are explicitly specified in the Appendix. 

\begin{figure}[h!]
    \centering
    \includegraphics[width=1.0\linewidth]{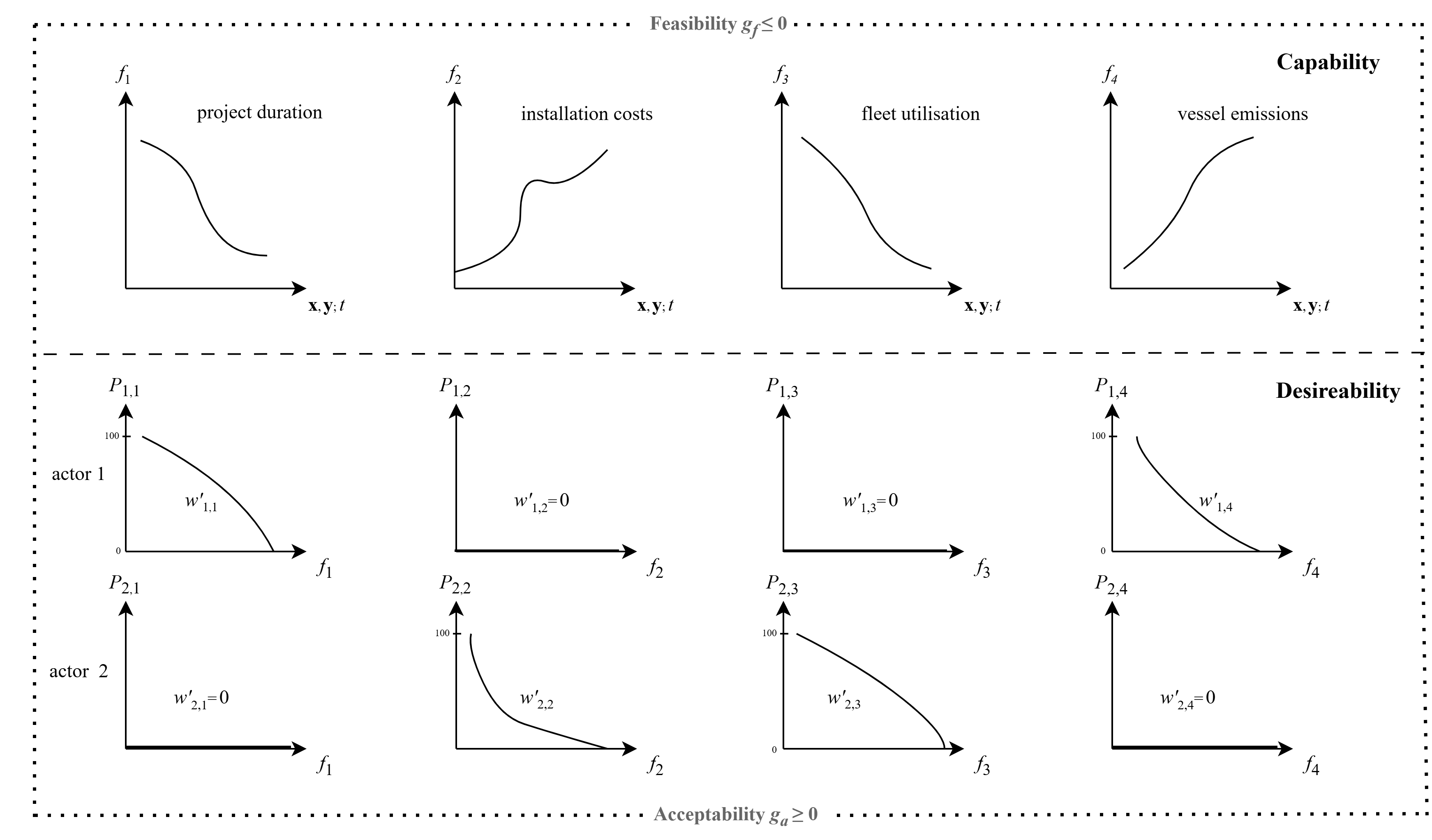}
     \vspace{-16pt}
    \caption{ODESYS–FIVES schematization (Example 1)}
    \label{fig:ODESYS_structureFIVES}
\end{figure}

\paragraph{(3) Desirability}
There are two actors in the system, each assigned an equal global weight \(w_k = 0.5\). Both actors may, in principle, express preferences over all four performance dimensions \(f_1\)–\(f_4\); however, in this example each focuses on the criteria most relevant to their operational interests. The energy service provider (\(k=1\)) assigns equal local weights to \(f_1\) and \(f_4\), i.e., \(w_{1,1} = w_{1,4} = 0.5\), while all other local weights are zero. The marine contractor (\(k=2\)) assigns equal local weights to \(f_2\) and \(f_3\), i.e., \(w_{2,2} = w_{2,3} = 0.5\), with all remaining weights equal to zero. This results in the effective contribution weights \(w'_{k,i} = w_k \cdot w_{k,i}\), with \(w'_{1,1} = w'_{1,4} = w'_{2,2} = w'_{2,3} = 0.25\), and all remaining weights equal to zero. The corresponding individually constructed preference functions are in this case monotonic and are schematically illustrated in Figure~\ref{fig:ODESYS_structureFIVES}, which depicts the resulting preference structure (for the full functional forms, see \parencite{Wolfert2023OpenDesignSystems, VanHeukelum2024}). \textsc{Note:} This example is intentionally symmetric and assumes a single active preference per performance dimension. However, the ODESYS methodology also supports asymmetric preference configurations, enabling heterogeneous (potentially conflicting) stakeholder priorities and differentiated weighting structures within a participatory decision-making framework (see \parencite{TeuberWolfert2024ConfrontingConflicts} for details on conjoint-analysis-based preference elicitation).

\paragraph{(4) Acceptability}
In this case, the minimum acceptable preference threshold in preference space is zero, i.e., \( P_{k,i} \in [0,100] \). Moreover, the actors do not introduce implicit (negotiable) acceptability limits through performance reference intervals, i.e., \( (f_i^{\mathrm{loc}},\, f_i^{\mathrm{upc}}) = (f_i^{\min},\, f_i^{\max}) \). Consequently, the actor-defined acceptability conditions do not restrict the admissible design--decision space, resulting in \( \mathcal{S}_{f,a} = \mathcal{S}_f \). The acceptability solution space (for \(K = 2\) and \(I = 4\)), prior to optimisation over the design–decision vector \( \mathbf{x} \), is formally defined as:
\[
\mathcal{S}_a
=
\left\{
\mathbf{x}
\;\middle|\;
g_a^{(n)} \ge 0,\;
f_i(\mathbf{x}) \in [f_i^{\min},\,f_i^{\max}],
\; n = 1..8
\right\}.
\]

Accordingly, all preference functions satisfy \( P_{k,i} \in [0,100] \), with 100 and 0 corresponding to the relatively best and worst feasible performance values, respectively. In this case, the preference-function endpoints coincide with the feasible performance extrema, i.e., \( P_{k,i} = 100 \) at \( f_i = f_i^{\min} \) (or \( f_i^{\max} \), depending on preference direction) and \( P_{k,i} = 0 \) at the opposite bound. Hence, all feasible system performances are admissible and fully represented in the preference space (see the schematic representation in Figure~\ref{fig:ODESYS_structureFIVES}; for further details, see \parencite{Wolfert2023OpenDesignSystems, VanHeukelum2024}).
\textsc{Note:} In general, actors may introduce explicit or implicit acceptability constraints—such as limits on project duration, costs, or emissions—which restrict the admissible range of the corresponding performance functions. In such cases, the acceptability space becomes a proper subset of the feasibility space, reflecting performance levels that are feasible but no longer considered acceptable.

\paragraph{(5) Solvability}
Given the defined acceptability–feasibility and preference–performance structure, the optimal design–decision vector \( \mathbf{x}^* \) can be obtained by solving the IMAP–MODO preference-based multi-objective design optimisation problem for multi-criteria group decision-making (see Eq.~\eqref{maxZ}). In this case, \(\mathcal{S}_{f,a}\) is bounded solely by the feasibility constraints \(g_f\), while the preference functions span the full range \([0,100]\). The objective is to identify the best-fit-for-common-purpose solution that maximises the aggregated preference score. This optimisation can be performed using an appropriate instantiation from the Inter-Generational Solver (IGS) family; see the Data Availability Statement for implementation details. As this paper focuses on methodology, no numerical results are presented here. Detailed applications and results for this problem setting can be found in \parencite{Wolfert2023OpenDesignSystems, VanHeukelum2024, VanHeukelum2023, Teuber2025Odycon}. These studies show that IMAP yields superior design–decision outcomes and, at minimum, matches established decision-valid methods for group decision-making (e.g., min–max), while outperforming single-objective optimisation approaches. Moreover, \parencite{Ozer2026} demonstrates that, for the considered demonstrators of \parencite{Wolfert2023OpenDesignSystems} (including this case), the Pareto front primarily constitutes an ordering of solutions rather than a decision instrument, as the best-fit-for-common-purpose solution does not generally belong to the Pareto set in the presence of non-monotonic real-world preferences or conflicting actor preferences for a single performance dimension.

\subsubsection*{\textsc{Demonstrator}~\#2: Multi-Constrained Vessel Allocation Problem}
\vspace*{-0.6em}
Within Boskalis, the ODESYS methodology has been applied to complex operational design--decision problems, including a highly constrained dynamic vessel allocation setting. This led to the development of the AlloDyn decision-support tool, embedded in the Vessel Operations Management System (VOMS) of the Marine Services business unit, which manages a global fleet of approximately 30 vessels and allocates incoming commercial projects to operationally available and capable vessels. This integrative decision problem addresses vessel selection and scheduling under multiple operational constraints, optimising sailing distance, costs, fuel consumption, and make-span. Two participating actors represent the operational and commercial teams with potentially conflicting objectives: the Operations team seeks to minimise operational costs, including mobilisation distance and fuel consumption, while the Commercial team aims predominantly to maximise project throughput and thereby minimise overall portfolio make-span. This results in four system performance functions in total. The industrial AlloDyn decision problem can be formulated as an ODESYS operator in Equation~\eqref{OD operator}, see the Appendix. While the full industrial problem considers \(k=1..2\) actors and \(i=1..4\) performance functions, the present demonstrator focuses on the methodological aspects of the MODO ODESYS/FIVES framework rather than on the industrial optimisation results. These results are reported in \parencite{Timp2026} through a reduced two-dimensional performance-space representation of the HVASP allocation problem. Accordingly, the performance space is here also reduced and re-indexed to the two dominant actor-specific performance functions, namely total cost (\(f_1\)) for the Operations BU and total make-span (\(f_2\)) for the Commercial BU, yielding \(i=1..2\).  In a quasi-dynamic setting (e.g.\ where each week a new rolling vessel allocation is required), time represents the evaluation instant $T^*$ of a fixed ODESYS system state and is therefore omitted, yielding $\mathbf{F}(\mathbf{x},\mathbf{y})$ evaluated at $t = T^*$. \vspace*{0.6em}

The ODESYS structure in FIVES is defined after the Note below.\vspace*{0.6em}

\textsc{Note:} This formulation represents a highly constrained combinatorial decision problem, characterised by tightly coupled sequencing, activity, and resource constraints. In practice, such problems are typically addressed using constraint programming (CP) or mixed-integer (non-linear) programming (MINLP) techniques. In this work, however, the problem is embedded within the ODESYS/FIVES framework to enable preference-based group decision-making across multiple actors and non-linear performance dimensions, rather than purely optimisation-driven formulations. Conventional CP, MIP, and MINLP approaches primarily focus on numerical optimisation and do not explicitly model actor preferences as a governing element of the decision problem, which might finally lead to invalid decision support (see Introduction).

\paragraph{(1) Capability}
The system capability is described by a two-dimensional performance vector of interest (as in the HVASP problem)
\[
\mathbf{F}(\mathbf{x},\mathbf{y}) =
\big[
f_1(\mathbf{x},\mathbf{y}),
f_2(\mathbf{x},\mathbf{y})
\big],
\]
representing total cost \(f_1(\mathbf{x},\mathbf{y})\) (including mobilisation and standby costs required to complete the activity portfolio), and total make-span \(f_2(\mathbf{x},\mathbf{y})\) (i.e.\ the total time required to complete the full set of activities). The controllable design–decision vector is \(\mathbf{x} = (x_1,\ldots,x_6)\), while the uncontrollable parameter vector is \(\mathbf{y} = (y_1,\ldots,y_{15})\). The full analytical expressions for cost \(f_1\) and make-span \(f_2\), including their dependence on \(\mathbf{x}\) and \(\mathbf{y}\), as well as the domain constraints \(g_f^{(0)}(\mathbf{x},\mathbf{y})\), are provided in the Appendix (note: in the full AlloDyn formulation, \(f_2\) corresponds to cost and \(f_4\) to make-span).

\paragraph{(2) Feasibility}
System feasibility is enforced through a set of constraints that ensure:
- domain validity (variable bounds on the decision vector $\mathbf{x}$ and exogenous parameter vector $\mathbf{y}$),
- activity validity (feasible allocation of vessels to activities),
- sequencing consistency (activities follow valid temporal orders),
- and path continuity (feasible travel paths between activity locations). 
The feasible system solution space prior to optimisation over the design–decision vector $\mathbf{x}$ is defined as:

\[
\label{eq:system_feasibility}
\mathcal{S}_f
=
\left\{\mathbf{x}
\;\middle|\;
g_f^{(m')} \le 0,
\; m' = 0..13
\right\}.
\]

Here, $g_f^{(0)}$ represents the domain constraints for the vectors $\mathbf{x}$ and $\mathbf{y}$, while constraints $g_f^{(1)}$–$g_f^{(13)}$ capture activity, sequencing, and path constraints. All individual constraints are explicitly specified in the Appendix.

\paragraph{(3) Desirability}
There are two actors in the system, each assigned an equal global weight \( w_k = 0.5 \). Both actors may, in principle, express preferences over the two performance dimensions \( f_1 \) (cost) and \( f_2 \) (make-span); however, in this demonstrator each prioritises these dimensions differently. The operations team (\( k = 1 \)) assigns a dominant local weight to cost, i.e., \( w_{1,1} = 0.9 \), and a minor weight to make-span, i.e., \( w_{1,2} = 0.1 \). The commercial team (\( k = 2 \)) assigns a dominant local weight to make-span, i.e., \( w_{2,2} = 0.8 \), and a secondary weight to cost, i.e., \( w_{2,1} = 0.2 \). This results in the effective contribution weights \( w'_{k,i} = w_k \cdot w_{k,i} \), yielding \( w'_{1,1} = 0.45 \), \( w'_{1,2} = 0.05 \), \( w'_{2,2} = 0.4 \), and \( w'_{2,1} = 0.1 \). The corresponding preference functions are in this case assumed to be linear.

A key characteristic of this demonstrator is that the actors exhibit conflicting preference directions over the make-span performance dimension: the operations team prefers longer durations due to perceived execution safety, whereas the commercial team prefers shorter durations due to higher revenue potential. Both actors, however, share a common preference orientation for cost, namely that lower costs are preferred, although they assign different levels of importance to this criterion. These opposing (partly aligned and partly conflicting) linear preference structures reflect structurally different operational and commercial interests. In addition, a coupling between performance dimensions is present: reducing the make-span (e.g., through faster sailing) leads to increased costs, further reinforcing the trade-off between the two criteria. \textsc{Note:} This fully individualised form of participation in MODO--IMAP, both in terms of stakeholder interests and weighting, represents a novel feature compared to often non-participatory classical multi-objective optimisation (MOO) approaches.

\textsc{Note:} This fully individualised form of participation in MODO--IMAP, both in terms of stakeholder interests and weighting, represents a novel feature compared to often non-participatory classical multi-objective optimisation (MOO) approaches. The preference structure in this demonstrator is intentionally simplified to highlight actor-specific priorities and conflicting preferences using linear functions. In general, ODESYS supports asymmetric and non-linear preference configurations, enabling heterogeneous stakeholder interests and differentiated weighting structures within a participatory decision-making framework (see \parencite{TeuberWolfert2024ConfrontingConflicts, Teuber2025Odycon}).

\paragraph{(4) Acceptability}
In this case, the minimum acceptable preference threshold in preference space is zero, i.e., \( P_{k,i} \in [0,100] \). Moreover, the actors introduce implicit (negotiable) acceptability limits through their choice of performance reference intervals. In particular, the operations team imposes an upper bound on cost, \( f_1^{\mathrm{upc}} = \overline{B} < f_1^{\max} \), while the commercial team imposes an upper bound on make-span, \( f_2^{\mathrm{upc}} = \overline{M} < f_2^{\max} \). Consequently, the actor-defined acceptability conditions restrict the admissible design--decision space, resulting in \( \mathcal{S}_{f,a} \subset \mathcal{S}_f \). The acceptability solution space (for \(K = 2\) and \(I = 2\)), prior to optimisation over the design–decision vector \( \mathbf{x} \), is formally defined as:
\[
\mathcal{S}_a
=
\left\{
\mathbf{x}
\;\middle|\;
g_a^{(n')} \ge 0,\;
f_i(\mathbf{x}) \in [f_i^{\mathrm{loc}},\,f_i^{\mathrm{upc}}],
\; n' = 1..4
\right\}.
\]

Accordingly, all preference functions satisfy \( P_{k,i} \in [0,100] \), with 100 and 0 corresponding to the relatively best and worst acceptable performance values, respectively. In particular, \( P_{1,1} = 0 \) at \( f_1 = \overline{B} \) and \( P_{2,2} = 0 \) at \( f_2 = \overline{M} \), while the remaining preference functions attain their endpoints at the corresponding feasible performance extrema \( f_i^{\min} \) and \( f_i^{\max} \). Hence, only performance values within these actor-defined intervals are admissible, and the preference space represents a restricted subset of the feasible system performances (for further detailed results, see the HVASP problem in \parencite{Timp2026}).

\paragraph{(5) Solvability}
Given the defined acceptability–feasibility and preference–performance structure, the optimal design–decision vector \( \mathbf{x}^* \) can be obtained by solving the IMAP–MODO preference-based multi-objective design optimisation problem for multi-criteria group decision-making (see Eq.~\eqref{maxZ}). In this case, \(\mathcal{S}_{f,a}\) is restricted by both feasibility constraints \(g_f\) and actor-defined acceptability conditions \(g_a\), while the preference functions span the admissible interval \([0,100]\) within the defined reference domains. The objective is to identify the best-fit-for-common-purpose solution that maximises the aggregated preference score.

This optimisation can be performed using an appropriate instantiation from the Inter-Generational Solver (IGS) family. For this demonstrator, an intergenerational Biased Random-Key Genetic Algorithm (IGS-BRKGA) is employed to effectively explore the highly constrained solution space. The algorithm adopts a decoder-based approach to construct feasible schedules and incorporates constraint handling through preference-consistent evaluation, enabling efficient identification of admissible and high-quality solutions within \( \mathcal{S}_{f,a} \).

As this paper focuses on methodology, no numerical results are presented here. Detailed applications and results for this problem setting can be found in the HVASP problem of \parencite{Timp2026}. These studies confirm that IMAP yields robust design–decision outcomes by directly operating in the integrated performance–preference space, rather than relying on classical objective-space trade-offs that may yield numerical optima but do not necessarily produce design–decision-valid outcomes.

Finally, a comparative performance evaluation of IMAP–IGS lies beyond the scope of this paper. In particular, a formal assessment on standard benchmark MOO problems—such as (i) the DAS-CMOP constrained multi-objective benchmark suite, and (ii) the MO-VRPTW multi-objective vehicle routing problem with time windows—is not included here. This evaluation, including a classification of design–decision problem characteristics under which a posteriori Pareto analysis can yield meaningful results, has recently been completed and will be reported in forthcoming work \parencite{Timp2026}. Prior work \parencite{Wolfert2023OpenDesignSystems} indicates that IMAP and the IGS-GA family achieve competitive or superior design–decision outcomes compared to established decision-valid methods for group decision-making (e.g., min–max approaches), while performing at least as well as single-objective and a posteriori Pareto-set-based approaches.

\subsection*{4. Conclusions}

The demonstrative examples above illustrate that conventional multi-objective approaches, when applied without integrated preference modelling, cannot guarantee a unique, consistent, and mathematically valid group decision outcome. In contrast, the ODESYS structure in FIVES implements a pure performance–preference paradigm, where (1) capability (system performance), (2) feasibility (hard constraints), (3)  desirability (stakeholder preferences) and (4) acceptability (negotiable constraints), are fully integrated within a single associative decision space, while (5) solvability is achieved by integratively maximising the aggregated preference (IMAP). For this, preferences are mapped into a single preference domain and aggregated via a uniquely defined, PFM-consistent linear aggregation operator, ensuring that all relative differences are preserved and that the resulting group-optimal solution simultaneously satisfies the four MOO–MCDM conditions—Preference-Key, Integration, Association, and Uniqueness—yielding a single, consistent, and fully group decision-valid outcome. By contrast, all other methods reviewed in the references cited in this paper—including CP-, MIP-, and conventional MOO approaches—do not satisfy the four MOO–MCDM conditions and therefore do not guarantee a decision-valid solution. This limitation is consistent with broader observations in the literature regarding the lack of unified and methodologically consistent MOO–MCDM integration \parencite{Pajasmaa2025GroupDecision, Ferdous2024, french2023reflections}. Moreover, many decision models rely on mathematical constructs that do not satisfy measurement-theoretic requirements, and therefore may produce misleading optimal solutions. ODESYS addresses and resolves these limitations through its integrated preference-based framework. Importantly, the framework preserves individual design freedom within the feasible and acceptable solution space, enabling actors to express their priorities while contributing to a collective best-fit, resulting in a synthesis of common purpose across the system’s 'infinite' physical possibilities. 

Through its structure-in-FIVES formulation, ODESYS operationalises this methodology for complex, multi-system design and decision problems, rigorously enforcing both hard feasibility constraints and soft acceptability conditions. IMAP identifies the unique best-fit-for-common-purpose design–decision vector using an intergenerational genetic search algorithm (GA), accommodating multiple actors with diverse preferences. In doing so, it finds the outcome with maximum aggregated preference, maintains full traceability and participatory engagement and enables forward-looking design decision-making rather than selecting from a curated numeric set, which unveils the illusion of free choice in conventional MCDM. 

Taken together, these results establish the IMAP–IGS solver and the ODESYS/FIVES framework as a rigorous, preference-consistent, and operationally validated methodology for multi-objective, multi-actor design–decision-making. Central to this approach is the principle that one can only find a design–decision solution once its capability–desirability has been defined, guiding the search toward the maximisation of integrated aggregated preference within a feasible and acceptable design-decision space, rather than relying on numerical ordering within objective space. The framework demonstrates how system performance (\textit{capability}), stakeholder preference (\textit{desirability}), and the distinction between \textit{feasibility} and \textit{acceptability} can be coherently integrated to produce decisions that are simultaneously transparent, associative, and \textit{solvable} as a best-fit-for-common-purpose. Within ODESYS/FIVES, IMAP performs comparably to or better than both classical multi-objective and single-objective optimisation approaches, consistently yielding a single, preference-consistent design–decision outcome \parencite{Wolfert2023OpenDesignSystems, VanHeukelum2024}. Consequently, when individual actors are willing to relax their purely self-interested positions, the group outcome achieves a higher aggregated preference, transforming creative conflicts into a best-fit synthesis—a shared ‘yes’—without compromising the physical realities of the system. This positions ODESYS as a framework capable of structuring and confronting complex socio-technical design–decision problems in a principled and operational manner.

In this sense, ODESYS/FIVES reorients MOO–MCDM paradigms towards the integrative maximisation of aggregated preference. After all, value is about more than money. Rather than searching for (monetised) objective-space optimality without decision validity, ODESYS aims for a best-fitting synthesis between performance possibilities and preference-based purpose, within the boundaries of feasibility and acceptability in the design–decision space. By leveraging the synergy of ‘systems-thinking-social’ and ‘design-thinking-slow’, it delivers a best-fit-for-and-from-the-whole within reach: a conscience of freedom!\\

\section*{Disclosure Statement}
The author declares that there are no competing interests.

\section*{Acknowledgements}

The author would like to thank his former TUDelft students and colleagues at Boskalis, and in particular Harold van Heukelum, Lukas Teuber, Emre \"Ozer, and Lennard Timp. Without their open design learning mindset and substantive contributions, the ODESYS framework would not have reached its current level of maturity.

Finally, the author gratefully acknowledges all users of a wide range of ODESYS-based applications at Boskalis and other industrial organisations, as well as the MSc students at TU Delft, for their indirect validation of the ODESYS approach over many years. This continued deployment in industry and academia provides a robust and continuously expanding real-world application and validation basis for ODESYS, for which the author expresses his sincere gratitude.

\section*{Data Availability}

The open-source Preferendus software for preference-based optimisation used for the IMAP demonstrators in the ODESYS book \parencite{Wolfert2023OpenDesignSystems} is publicly available via the GitHub repository:
\url{https://github.com/TUDelft-Odesys}, including the early IMAP/IGS 1.0 implementation. The results of the first demonstrator presented in this paper (the floating wind farm demonstrator) are described in Chapter~8.5 of the ODESYS book and can also be accessed through the same GitHub repository. Furthermore, extended variants of this design–decision example are available in \parencite{Teuber2025Odycon}, \parencite{VanHeukelum2024} or \parencite{VanHeukelum2023}.

A contemporary and further extended implementation of Preferendus, incorporating a collection of Inter-Generational Solver (IGS) instantiations with problem-specific decoders, is available at:
\url{https://github.com/TUDelft-IMAP-IGS}. 

Within this IGS family, IGS-BRKGA is used for the HVASP and MO-VRPTW problems, each with dedicated decoder structures, and IGS-GA-II is applied to the DAS-CMOP problem (see  \parencite{Timp2026}). In addition, IGS-SI is used for a routing problem presented in \parencite{Ozer2026}. 

In particular, the IGS-BRKGA instantiation is used for the second demonstrator in this paper (the dynamic vessel allocation problem). The corresponding industrial Allodyn results are subject to confidentiality constraints and are available upon reasonable request. However, an academic variant of this highly constrained combinatorial design--decision problem is available in \parencite{Timp2026}, corresponding to the HVASP problem.

\newpage
\section*{Appendix}

\setcounter{equation}{0}
\renewcommand{\theequation}{A\arabic{equation}}
 
Let an integrative design and group decision problem for a given system state be defined by the ODESYS system operator as the basis for it's structure in FIVES as:
 
\begin{equation}
\mathrm{OD}\Big(
P_{k,i}\big(\mathbf{F}(\mathbf{x},\mathbf{y})\big);\;
w'_{k,i}
\Big),
\end{equation}
 
Here, $\mathbf{F} = \big(f_1,\ldots,f_I\big)$ denotes the system performance (capability) vector, where each function $f_i(\mathbf{x},\mathbf{y})$ represents the $i$th performance dimension, and $I$ is the total number of performance functions. The functions $P_{k,i}(\cdot)$ define actor-specific preference functions over performance dimension $i$ for actor $k$, while $w'_{k,i}$ denote the corresponding local preference weights.
 
\subsection*{Example (1)}
This example illustrates the formulation of system capability functions for the floating wind farm integrative design--decision problem.

\subsubsection*{Capability}
For the offshore anchor installation problem, multi-system performance is captured by the following capability functions:

\begin{equation}
\mathbf{F} =
\big[
f_1(\mathbf{x}),
f_2(\mathbf{x},\mathbf{y}),
f_3(\mathbf{x}),
f_4(\mathbf{x},\mathbf{y})
\big].
\end{equation}

The project duration (total project make-span) is determined through a discrete-event simulation (DES):

\[
f_1(\mathbf{x}) = ft_{\text{proj}}(t_{1..3}(\mathbf{x}))
\]

where $ft_{\text{proj}}$ (in days) is obtained from the DES based on different vessel configurations $(x_1,x_2,x_3)$ and their operating times $t_{1..3}(\mathbf{x})$.

\vspace{6pt}

The total installation costs (EUR) are given by:
\[
f_2(\mathbf{x},\mathbf{y}) =
(815\,M_a + 40{,}000)\,y_6
+
\sum_{i=1}^{3} x_i\, t_i(\mathbf{x})\, y^{R}_i
\]

where the cost per anchor consists of a fixed component (€40,000 per anchor) and a variable component proportional to material usage (€815 per tonne), with \(M_a\) denoting the anchor mass (in tonnes), given by:

\[
M_a =
\left(
\pi x_5 x_4 + \frac{\pi}{4} x_4^2
\right) y_{7}
\]

\vspace{6pt}

Fleet utilisation is defined as:

\[
f_3(\mathbf{x}) =
\prod_{i=1}^{3} p_i^{x_i}
\]

\vspace{6pt}

The vessel emissions are:

\[
f_4(\mathbf{x},\mathbf{y}) =
\sum_{i=1}^{3} x_i\, t_i(\mathbf{x}) y^{R}_i\
\]

\vspace{12pt}

Here, the controllable decision vector is defined as $\mathbf{x} = (x_1,\dots,x_5)$, representing vessel allocation and anchor design variables.  
The uncontrollable parameter vector is defined as $\mathbf{y} = (y_1,\dots,y_10)$, representing environmental, geotechnical, and system-level design conditions assumed fixed per scenario. These include the loading conditions acting on the anchor as well as the physical characteristics of the soil and mooring system.

\begin{table}[h!]
\centering
\small
\caption{Domain constraints for the decision vector $\mathbf{x} = (x_1,\dots,x_5)$}
\begin{tabular}{c l l}
\toprule
\textbf{$\mathbf{x}$} & \textbf{Description} & \textbf{$g_f^{(0)}(x_i)$} \\
\midrule
$x_1$ & Number of small OCVs & $0 \le x_1 \le 3,\;\; x_1 \in \mathbb{Z}_{\ge 0}$ \\
$x_2$ & Number of large OCVs & $0 \le x_2 \le 2,\;\; x_2 \in \mathbb{Z}_{\ge 0}$ \\
$x_3$ & Number of barges & $0 \le x_3 \le 2,\;\; x_3 \in \mathbb{Z}_{\ge 0}$ \\
$x_4$ & Anchor diameter & $1.5 \le x_4 \le 4$ \\
$x_5$ & Anchor penetration length & $2 \le x_5 \le 8$ \\
\bottomrule
\end{tabular}
\label{tab:x_domain_example1}
\end{table}

\begin{table}[h!]
\centering
\small
\caption{Parameter vector $\mathbf{y} = (y_1,\dots,y_7,\mathbf{y}^{R},\mathbf{y}^{E})$ and domain constraints}
\begin{tabular}{c l l}
\toprule
\textbf{$\mathbf{y}$} & \textbf{Description} & \textbf{$g_f^{(0)}(y_i)$} \\
\midrule
$y_1$ & Working point force on anchor $F_a$ & $y_1 \ge 0$ \\
$y_2$ & Mooring configuration & $y_2 \in \text{Configurations}$ \\
$y_3$ & Anchor type & $y_3 \in \text{AnchorTypes}$ \\
$y_4$ & Soil conditions & $y_4 \in \text{SoilTypes}$ \\
$y_5$ & Mooring line properties & $y_5 \in \text{LineProperties}$ \\
$y_6$ & Number of anchors & $y_6 = 108$ \\
$y_7$ & Areal mass density of steel  & $y_7 \ge 0$ \\
$\mathbf{y}^{R}$ & Vessel day rates & $(y^{R}_1,y^{R}_2,y^{R}_3) = (R_1,R_2,R_3)$ \\
$\mathbf{y}^{E}$ & Vessel emission rates & $(y^{E}_1,y^{E}_2,y^{E}_3) = (E_1,E_2,E_3)$ \\
\bottomrule
\end{tabular}
\label{tab:y_domain_example1}
\end{table}

\vspace{-1.0em}
\paragraph{Notes:}
\begin{itemize}
    \item The soil is assumed to be clay with undrained shear strength $s_u = 60$ kPa and submerged unit weight $\gamma' = 9$ kN/m$^3$. The shaft--soil interface friction coefficient is $\alpha = 0.64$, and the seabed--chain friction coefficient is $\mu = 0.25$.
    \item The mooring line chain (diameter $240$ mm), attached at $0.5\,x_5$ penetration depth.
    \item The active bearing area coefficient is $A_{WB} = 2.5$, while $y_7$ denotes the areal mass density of steel, defined as $y_7 = \rho_{\text{steel}} \cdot h$ (approx.\ $78.5\,\text{kg/m}^2$).
    \item The DES provides vessel operating times $t_i(\mathbf{x})$, which are used in $f_1$, $f_2$, and $f_4$.
    \item Vessel day rates are $R_1 = 47{,}000$, $R_2 = 55{,}000$, and $R_3 = 35{,}000$ (EUR/day). Vessel emission rates are $E_1 = 30$, $E_2 = 40$, and $E_3 = 35$ (tonne/day).
    \item Fleet utilisation probabilities are $p_1 = 0.7$, $p_2 = 0.8$, and $p_3 = 0.5$.
\end{itemize}

\subsubsection*{Feasibility}
The feasibility system solution space prior to optimisation over the design–decision vector $\mathbf{x}$ is defined as:
\[
\mathcal{S}_f
=
\left\{
\mathbf{x}
\;\middle|\;
g_f^{(m)} \le 0,\;
m = 0,1,2
\right\},
\]

This formulation provides concrete instantiations of the feasibility operator. Here $g_f^{(0)}$ denotes the domain constraints (see Table~1 for their values). The logical constraints are given by
\[
g_f^{(1)}(\mathbf{x})
:=
-(x_1+x_2+x_3)+1
\le 0,
\]
ensuring the logical feasibility of the vessel-selection decision, and the physical constraints are given by
\[
g_f^{(2)}(\mathbf{x})
:=
F_a(y_{1..5}) - R_a(x_4,x_5,y_{1..5})
\le 0,
\]
ensuring the physical (load-bearing capacity) feasibility of the cable-anchor design.\\

\subsection*{Example (2)}
This example illustrates the formulation of system capability and feasibility functions for a dynamic vessel allocation problem with multiple constraints.
\subsubsection*{Capability}

For the vessel allocation problem, multi-system performance is captured by the following capability functions:
\begin{equation}
\mathbf{F} =
\big[
f_1(\mathbf{x},\mathbf{y}),
f_2(\mathbf{x},\mathbf{y}),
f_3(\mathbf{x},\mathbf{y}),
f_4(\mathbf{x},\mathbf{y})
\big].
\end{equation}
Mobilization distance is defined as:
\[
f_1(\mathbf{x},\mathbf{y}) =
\sum_{r \in R} \sum_{r' \in R}
x_4[r, r'] \cdot
y_{10}[\ell^{\text{end}}(r), \ell^{\text{start}}(r')]
\]

while the total cost (i.e., the sum of all mobilization and standby costs, excluding the one-to-one activity sailing costs, to complete the activity portfolio) is given by:
\[
f_2(\mathbf{x},\mathbf{y}) =
\sum_{r \in R} \sum_{r' \in R}
 x_4[r, r']
\cdot \gamma \left(x_3[r], \ell^{\text{end}}(r), \ell^{\text{start}}(r'), x_6[r, r'], \delta_{r, r'} \right).
\]

Fuel consumption is modelled as:
\[
f_3(\mathbf{x}, \mathbf{y}) = 
\sum_{r \in R} \sum_{r' \in R} x_4[r, r']
\cdot y_{12}[x_3[r]](x_6[r, r']) \cdot \theta(x_3[r], \ell^{\text{end}}(r), \ell^{\text{start}}(r'), x_6[r, r'])
\]
and the total make-span (i.e., the total time required to complete the entire portfolio of activities) is expressed as:

\[
f_4(\mathbf{x},\mathbf{y}) =
\max_{a \in A}(x_1[a] + y_1[a]) - \min_{a \in A}(x_1[a])
\]

Here, the controllable decision vector is defined as $\mathbf{x} = (x_1,\dots,x_6)$, with components corresponding to timing, location, vessel assignment, sequencing, and sequence initiation decisions; the associated domain constraints are summarised in Table~\ref{tab:domain_constraints_simplified}. The uncontrollable parameter vector is defined as $\mathbf{y} = (y_1,\dots,y_{15})$, where each parameter captures exogenous system characteristics such as activity durations, time windows, spatial information, precedence relations, vessel capabilities, and sailing characteristics, and whose associated domains and feasibility constraints are listed in Table~\ref{tab:parameter_vector_constraints}. Finally, the index sets used in the formulation are defined in Table~\ref{tab:sets}, representing the vessels, activities, roles, and geographical locations relevant to the problem; the activity set is further partitioned by type, yielding a more compact role-based formulation than a binary vessel–activity assignment matrix.

\begin{table}[h!]
    \centering
    \small
    \caption{Domain constraints for the simplified decision vector $x = (x_1, \dots, x_6)$}
    \begin{tabular}{c l l}
    \toprule
    \textbf{$\mathbf{x}$} & \textbf{Description} & \textbf{$g_f^{(0)}(x_i)$} \\
    \midrule
    $x_1$ & Start time of activity $a \in A$ & $\underline{T}_a \le x_1[a] \le \overline{T}_a$ \\
    $x_2$ & Location choice for maintenance activity $a \in A_{maint}$ & $x_2[a] \in L_a$ \\
    $x_3$ & Vessel assigned to role $r \in R$ & $x_3[r] \in \mathcal{D}_r$ \\
    $x_4$ & Sequencing variable: role $r'$ follows $r$ & $x_4[r,r'] \in \{0,1\}$ \\
    $x_5$ & Sequence start indicator: role $r$ is first in sequence & $x_5[r] \in \{0,1\}$ \\
    $x_6$ & Average speed for sub-route $(r, r')$ & $x_6[r, r'] \in [\underline{s}, \overline{s}]$\\
    \bottomrule
    \end{tabular}
    \label{tab:domain_constraints_simplified}
\end{table}

\begin{table}[H]
    \centering
    \small
    \caption{Parameter vector $\mathbf{y} = (y_1, \dots, y_{15})$ and their domain constraints}
    \begin{tabularx}{\textwidth}{c l X}
    \toprule
    \textbf{$\mathbf{y}$} & \textbf{Description} & \textbf{$g_f^{(0)}(y_i)$} \\
    \midrule
    $y_1$ & Duration of activity $a$ & \(y_1[a] \ge 0\) (Days) \\
    $y_2$ & Start time window for activity $a$ & \(y_2[a] = [\underline{y_2[a]}, \overline{y_2[a]}] \ge 0\) \\
    $y_3$ & Start location for towing activity $a \in A_{tow}$ & \(y_3[a] \in \text{Locations}\) \\
    $y_4$ & End location for towing activity $a \in A_{tow}$ & \(y_4[a] \in \text{Locations}\) \\
    $y_5$ & Locations for maintenance activity $a \in A_{maint}$ & \(y_5[a] \subseteq \text{Locations}\) \\
    $y_6$ & Predecessor of activity $a$ & \(y_6[a] \in A \cup \{\emptyset\}\) \\
    $y_7$ & Parent activity of role $r$ & \(y_7[r] \in A\) \\
    $y_8$ & Vessel domain for role $r$ & \(y_8[r] \subseteq V\) \\
    $y_9$ & Set of roles belonging to activity $a$ & \(y_9[a] = \{r \in R \mid \alpha(r)=a\}\) \\
    $y_{10}$ & Sailing distance from location $\ell$ to $\ell'$ & \(y_{10}[\ell,\ell'] \ge 0\) \\
    $y_{11}$ & Daily mobilisation rate for vessel $v$ & \(y_{11}[v] \ge 0\) \\
    $y_{12}$ & Fuel consumption rate function for vessel $v$ &
    \parbox[t]{0.75\textwidth}{%
    \(\displaystyle y_{12}[v](s) \ge 0,\;\forall s \in [y_{13}[v]^{\min}, y_{13}[v]^{\max}]\)
    } \\
    $y_{13}$ & Feasible sailing speed range for vessel $v$ &
    \parbox[t]{0.75\textwidth}{%
    \(\displaystyle y_{13}[v] = (s_v^{(0)}, s_v^{(1)}, \dots, s_v^{(|S_v|-1)}),\; s_v^{(k)} > 0\)
    } \\
    $y_{14}$ & Fuel price for vessel $v$ & \(y_{14}[v] \ge 0\) \\
    $y_{15}$ & Standby cost discount factor for vessel $v$ & \(y_{15}[v] \in [0,1]\) \\
    \bottomrule
    \end{tabularx}
    \label{tab:parameter_vector_constraints}
\end{table}

\begin{table}[H]
    \centering
    \small
    \caption{All sets in the fleet allocation problem.}
    \begin{tabular}{l l}
    \toprule
    \textbf{Set} & \textbf{Description} \\
    \midrule
    $V = \{v_0, v_1, \ldots, v_n\}$ & Set of $n$ vessels \\
    $A = \{a_0, a_1, \ldots, a_m\}$ & Set of $m$ activities \\
    $A^{\text{tow}} \subseteq A$ & Subset of towing activities \\
    $A^{\text{maint}} \subseteq A$ & Subset of maintenance activities \\
    $R = \{r_0, r_1, \ldots, r_p\}$ & Set of $p$ roles \\
    $L = \{\ell_0, \ell_1, \ldots, \ell_q\}$ & Set of $q$ locations \\
    \bottomrule
    \end{tabular}
    \label{tab:sets}
\end{table}
 
\paragraph{Notes:}
\begin{itemize}
    \item $\theta(v, \ell, \ell', s)$ is the sailing duration for vessel $v$ to travel from location $\ell$ to $\ell'$ at average speed $s$, defined by
    \[
    \theta(v, \ell, \ell', s) = \left\lceil \frac{y_{10}[\ell, \ell']}{24 \cdot s} \right\rceil
    \]
    \item $\delta_{r,r'} = x_1[y_7[r']]-\left( x_1[y_7[r]] + y_1[y_7[r]] \right)$ denotes the total transition time between role $r$ and $r'$, where $y_7[r]$ denotes the parent activity of role $r$.
    \item $\gamma(v, r, r', s, \delta_{r,r'})$ is the cost for vessel $v$ to travel from role $r$ to $r'$ at speed $s$ with a total travel duration of $\delta_{r, r'}$, defined by 
    \begin{align*}
\gamma(v, r, r', s, \delta_{r,r'}) &= 
\theta(v, \ell^{\text{end}}(r), \ell^{\text{start}}(r'), s) 
\cdot \Big(y_{11}[v] + y_{14}[v] \cdot y_{12}[v](s) - y_{15}[v] \cdot y_{11}[v] \Big) \\
&\quad + y_{15}[v] \cdot y_{11}[v] \cdot \delta_{r,r'}
\end{align*}
    \item $\underline{s} = \min_{v \in V} (y_{13}[v])$, $\overline{s} = \max_{v \in V}(y_{13}[v])$ are the min and max available speeds across all vessels.
    \item The start and end locations of role $r$ are defined as
    \[
    \ell^{\text{start}}(r) =
    \begin{cases}
    y_3[y_7[r]] & \text{if } y_7[r] \in A^{\text{tow}} \\[2mm]
    x_2[y_7[r]] & \text{if } y_7[r] \in A^{\text{maint}}
    \end{cases}
    \]
    \[
    \ell^{\text{end}}(r) =
    \begin{cases}
    y_4[y_7[r]] & \text{if } y_7[r] \in A^{\text{tow}} \\[2mm]
    x_2[y_7[r]] & \text{if } y_7[r] \in A^{\text{maint}}
    \end{cases}
    \]
    so that role locations are determined either by fixed parameters ($y_3, y_4$) for towing activities or by the decision variable $x_2$ for maintenance activities.
\end{itemize}
 
\subsubsection*{Feasibility }
The feasibility system solution space prior to optimisation over the design–decision vector $\mathbf{x}$ is defined as:
\begin{equation}
\mathcal{S}_{f}
=
\left\{
\mathbf{x}
\;\middle|\;
g_f^{(m')} \le 0,\;
m' = 0..13
\right\}.
\end{equation}

The feasibility operator $g_f^{(m')}$, representing the set of different feasibility constraints indexed by $m'$, is decomposed into the following constraint instantiations:
 
\paragraph {(a) activity constraints} No vessel may be assigned to more than one role within the same activity:
\[
g_f^{(1)}(\mathbf{x},\mathbf{y})
= x_3[r] - x_3[r']
\neq 0,
\quad
\forall r,r' \in R,\; r \neq r'.
\]
 
An activity cannot start before its predecessor has finished:
\[
g_f^{(2)}(\mathbf{x},\mathbf{y})
= x_1[a] - \big(x_1[y_6[a]] + y_1[y_6[a]]\big)
\ge 0,
\quad
\forall a \in A:\; y_6[a] \neq \emptyset .
\]
 
\paragraph {(b) sequencing constraints} In the sequence of activities there can be at most one successor:
\[
g_f^{(3)}(\mathbf{x},\mathbf{y})
=
\sum_{r' \in R} x_4[r,r'] - 1
\le 0,
\quad \forall r \in R,
\]
 
and at most one predecessor
\[
g_f^{(4)}(\mathbf{x},\mathbf{y})
=
\sum_{r \in R} x_4[r,r'] - 1
\le 0,
\quad \forall r' \in R.
\]
 
There can not be self-loops:
\[
g_f^{(5)}(\mathbf{x},\mathbf{y})
=
x_4[r, r]
\le 0,
\quad \forall r \in R.
\]

If role r' follows r, they must be assigned to the same vessel:

\[
g_f^{(6)} = x_4[r, r'] =1 \rightarrow  x_3[r] = x_3[r']
\]
 
Transitions between roles of the same activity are prevented, as such roles are
performed simultaneously rather than sequentially:
\[
g_f^{(7)}(\mathbf{x},\mathbf{y})
=
x_4[r,r']
\le 0,
\quad
\forall r,r' \in R:\; y_7[r] = y_7[r'].
\]
 
There can not be temporal precedence in sequences
\[
g_f^{(8)}(\mathbf{x},\mathbf{y})
=
x_1[y_7[r]] - x_1[y_7[r']]
< 0,
\quad \forall r,r' \in R:\; x_4[r,r'] = 1.
\]
 
There should be sufficient travel-time for consecutive roles:
\[
g_f^{(9)}(\mathbf{x},\mathbf{y})
=
x_1[y_7[r']] - 
\Big(
x_1[y_7[r]] + y_1[y_7[r]]
+ \theta(x_3[r], \ell^{\text{end}}(r), \ell^{\text{start}}(r'), y_{13}[x_3[r]]^{max})
\Big)
\le 0,
\]

where $y_{13}[x_3[r]]^{max}$ is the maximum feasible speed for the selected vessel, and for all \(r,r' \in R\) such that
\(x_4[r, r'] = 1\) and \(x_3[r] = x_3[r']\).
 
\paragraph {(c) path constraints}
 
Each vessel must have exactly one sequence start:
\[
g_f^{(10)}(\mathbf{x},\mathbf{y})
=
\sum_{r \in R}
[x_3[r] = v]\,
\left[\sum_{r' \in R} x_4[r',r] = 0\right]
= 1,
\quad \forall v \in V
\]
 
and exactly one sequence end:
\[
g_f^{(11)}(\mathbf{x},\mathbf{y})
=
\sum_{r \in R}
[x_3[r] = v]\,
\left[\sum_{r' \in R} x_4[r,r'] = 0\right]
= 1,
\quad \forall v \in V.
\]
 
For each vessel, the assigned roles must together form exactly one continuous path (path cardinality condition):
\[
g_f^{(12)}(\mathbf{x},\mathbf{y})
=
\sum_{r \in R} \sum_{r' \in R}
x_4[r,r'] \cdot [x_3[r] = v]
-
\max\!\left(
0,\;
\sum_{r \in R}[x_3[r] = v] - 1
\right)
= 0,
\quad \forall v \in V.
\]

If role $r'$ follows $r$ for vessel $v$, the chosen speed on that sub-route must be in the set of available sailing speeds of that vessel.
\[
g_f^{(13)}(\mathbf{x}, \mathbf{y}) = 
(x_4[r,r'] = 1 \;\land\; x_3[r] = v) \;\Rightarrow\; x_6[r,r'] \in y_{13}[v] \quad \forall v \in V,\; \forall r, r' \in R 
\]

\end{document}